\documentclass[moor]{informs1_for_arxiv}

\usepackage[sort]{natbib}
 \NatBibNumeric
 \def\bibfont{\small}%
 \bibpunct[, ]{[}{]}{,}{n}{}{,}%

\usepackage[colorlinks=true,breaklinks=true,bookmarks=true,urlcolor=blue,
     citecolor=blue,linkcolor=blue,bookmarksopen=false,draft=false]{hyperref}

\def\EMAIL#1{\href{mailto:#1}{#1}}
\def\URL#1{\href{#1}{#1}}         

\TheoremsNumberedThrough     

\EquationsNumberedThrough    

\usepackage[utf8]{inputenc}
\usepackage[T1]{fontenc} 
\usepackage{mathtools}
\usepackage{graphicx,subcaption}
\usepackage{amssymb,amsmath}
\usepackage{algorithm}
\usepackage{soul}
\usepackage{algorithmic}
\usepackage{epstopdf}
\usepackage{personal}
\usepackage{color}
\usepackage{answers}
\usepackage{setspace}
\usepackage{graphicx}
\usepackage{enumitem}
\usepackage{multicol}
\usepackage{mathrsfs}
\usepackage[margin=1in]{geometry} 
\usepackage{amsmath,amsfonts,amssymb,amsbsy,mathrsfs,bbm}
\usepackage{mathtools}
\usepackage{authblk}
\usepackage{breqn}
\usepackage{scalerel}
\usepackage[below]{placeins}

\newcommand{\inP}{\mbox{$\,\stackrel{\scriptsize{\mbox{p}}}{\rightarrow}\,$}}
\newcommand{\as}{\mbox{$\,\stackrel{\scriptsize{\mbox{\scriptsize wp1}}}{\rightarrow}\,$}}
\newcommand{\inD}{\mbox{$\,\stackrel{\scriptsize{\mbox{d}}}{\rightarrow}\,$}}

\newcommand{\Z}{\mathbb{Z}}

\newcommand{\E}{\mathbb{E}}

\renewcommand{\mathcal}{\mathscr}

\newcommand*\diff{\mathop{}\!\mathrm{d}}

\usepackage{tikz}
\usetikzlibrary{snakes}
\usetikzlibrary{arrows}
\usepackage[most]{tcolorbox}

\DeclarePairedDelimiterX{\expectarg}[1]{[}{]}{%
  \ifnum\currentgrouptype=16 \else\begingroup\fi
  \activatebar#1
  \ifnum\currentgrouptype=16 \else\endgroup\fi
}

\newcommand{\innermid}{\nonscript\;\delimsize\vert\nonscript\;}
\newcommand{\activatebar}{%
  \begingroup\lccode`\~=`\|
  \lowercase{\endgroup\let~}\innermid 
  \mathcode`|=\string"8000
}

\DeclarePairedDelimiter\ceil{\lceil}{\rceil}

\newcommand{\rv}[1]{{\color{blue}#1}}

\begin{document}


\ABSTRACT{Stochastic Gradient (SG) is the de facto iterative technique to solve stochastic optimization (SO) problems with a smooth (nonconvex) objective $f$ and a stochastic first-order oracle. SG's attractiveness is due in part to its simplicity of executing a \emph{single step} along the negative subsampled gradient direction to update the incumbent iterate. In this paper, we question SG's choice of executing a \emph{single step} as opposed to \emph{multiple steps} between subsample updates. Our investigation leads naturally to generalizing SG into Retrospective Approximation (RA) where, during each iteration, a ``deterministic solver'' executes possibly \emph{multiple steps} on a subsampled deterministic problem and stops when further solving is deemed unnecessary from the standpoint of statistical efficiency. RA thus formalizes what is appealing for implementation --- during each iteration, ``plug in'' a  solver, for example, L-BFGS line search or Newton-CG, \emph{as is}, and solve only to the extent necessary. We develop a complete theory using \emph{relative error} of the observed gradients as the principal object, demonstrating that almost sure and $L_1$ consistency of RA are preserved under especially weak conditions when sample sizes are increased at appropriate rates. We also characterize the iteration and oracle complexity (for linear and sub-linear solvers) of RA, and identify a practical termination criterion leading to optimal complexity rates. To subsume nonconvex $f$, we present a certain ``random central limit theorem'' that incorporates the effect of curvature across all first-order critical points, demonstrating that the asymptotic behavior is described by a certain mixture of normals. The message from our numerical experiments is that the ability of RA to incorporate existing second-order deterministic solvers in a strategic manner might be important from the standpoint of dispensing with hyper-parameter tuning.}

\KEYWORDS{Sample Selection, Machine Learning}
\MSCCLASS{}
\ORMSCLASS{Primary: ; secondary: }
\HISTORY{}

\RUNTITLE{}

\TITLE{A Retrospective Approximation Approach \\ for Smooth Stochastic Optimization}

\ARTICLEAUTHORS{
\AUTHOR{David Newton}
\AFF{Purdue University, Department of Statistics, \EMAIL{newton34@purdue.edu}, \URL{}}
\AUTHOR{Raghu Bollapragada}
\AFF{The University of Texas at Austin, Department of Mechanical Engineering, \EMAIL{Raghu.bollapragada@utexas.edu}, \URL{}}
\AUTHOR{Raghu Pasupathy}
\AFF{Purdue University, Department of Statistics; IIT Madras, Department of Comp. Sci. \& Engg., \EMAIL{pasupath@purdue.edu}, \URL{}}
\AUTHOR{Nung Kwan Yip}
\AFF{Purdue University, Department of Mathematics, \EMAIL{yipn@purdue.edu}, \URL{}}
} 




\newpage

\maketitle

\section{INTRODUCTION}\label{sec:intro}
We propose iterative algorithms to solve the unconstrained smooth stochastic optimization problem 
\begin{align*}\label{master}
 \underset{x \in \mathbb{R}^d}{\mbox{minimize }} \ & f(x) := \mathbb{E}[F(x,Y)] = \int F(x,y) \, P(dy),  \tag{$Q$}
\end{align*} 
where $Y$ is a random variable having distribution $P$ on a measurable space $(\mathcal{Y},\mathcal{A})$, and $F(\cdot,\cdot) : \mathbb{R}^d \times \mathcal{Y} \to \mathbb{R}$ is a ``random function''. Our standing assumption is that $F(x,Y)$ is a smooth function of $x\in\mathbb{R}^d$ for any $Y=y \in \mathcal{Y}$, 
and that we have access to a stochastic first-order blackbox oracle which allows for ``observing'' the function $F$ and its gradient $\nabla F$ at any $x \in \mathbb{R}^d$ and $Y=y \in \mathcal{Y}$. (We will make precise the notion of a stochastic first-order oracle and the underlying sigma algebra within the algorithm in Sections~\ref{sec:notation} and \ref{sec:ra} respectively.) Optimization problems of this form arise in a variety of applications arising within simulation optimization \cite{BCMS2018,Fu2005,Kim2014,Pasupathy2013,Pasupathy2018}, machine learning \cite{bottou2018optimization}, reinforcement learning \cite{Bertsimas2019b,Mania2018,Salimans2017}, data fitting \cite{Friedlander2012} and many others.

\subsection{Stochastic Gradient and Sample Average Approximation}
Stochastic Gradient (SG)~\cite{bottou2018optimization,2009judlannemsha} (originally called Stochastic Approximation in the seminal work of Robbins and Monro~\cite{1951robmon})  and Sample Average Approximation (SAA)~\cite{2009shadenrus,2015kimpashen} are the two main algorithmic frameworks used for solving \eqref{master}. (See~\cite{bottou2018optimization,2009judlannemsha,2015kimpashen,2009shadenrus} and references therein for extensive commentary on the respective strengths and limitations of SG and SAA.) In its most basic form, SG starts with an initial guess, and repeatedly updates the parameters $x$ by executing a gradient-descent like recursion using a sampled approximation of $f(x)$ at each iteration. Over the last decade or so, a number of variations of SG have been developed and demonstrated to be optimal. For example, mirror descent stochastic approximation~\cite{2009nemjudlansha} and accelerated stochastic approximation~\cite{2011lan} exhibit optimal complexity~\cite{2011lan} on stochastic convex optimization problems. Similarly, the $O(\epsilon^{-4})$ complexity exhibited by randomized stochastic gradient~\cite{2013ghalan} on nonconvex smooth stochastic optimization problems is optimal since it has been recently demonstrated~\cite{arjevani2019lower} that $O(\epsilon^{-4})$ forms a lower bound complexity for smooth nonconvex stochastic optimization using a first-order oracle.  

Despite such guarantees of optimality, work on SG has continued over the decades mostly in attempts to make SG work well in practice without the need for tuning hyperparameters. For example, the most recent wave of progress in SG has resulted in the incorporation of such features as momentum~\cite{2014kin}, adaptive stepsizes~\cite{2019lydfra,2018warwubot}, variance reduction techniques (SAG~\cite{2017schlerbac},
SAGA~\cite{2014defbaclac}, SVRG~\cite{2013johzha}), stochastic higher-order methods to incorporate curvature information via quasi-Newton update procedures \cite{Berahas2016,BollapragadaICML18,mokhtari2015global, schraudolph2007stochastic} or via subsampled Hessian approximations \cite{Bollapragada2018b,RoostaKhorasani2019,Pasupathy2018}, and adaptive sampling ideas~\cite{Bollapragada2018,Byrd2012,cartis2018global,Hashemi2014ada} aimed at gradually improving the quality of the approximation of $f$ as the iterates evolve. SG approaches, by their nature, tend to produce iterates that have low bias but high variance, low computational cost per iteration but inferior iteration complexity, and the need for tuning of hyper-parameters~\cite{2017cursch}. 

\begin{remark}
Note that we have refrained from using the popular misnomer ``Stochastic Gradient Descent'' in favor of the more appropriate ``Stochastic Gradient'' since SG is not a descent algorithm in that no guarantees of progress are available during each step.  
\end{remark}


A broad competing framework to SG is sample-average approximation (SAA). In contrast to SG, advances in SAA, especially its specialized optimization-oriented development, has taken place only since the 1990s~\cite{91heasch,90rubsha,shapiro4,1993sha,shapiro6,sha04,2009shadenrus,2001shakle} due largely to the seminal work of A. Shapiro. (We note that the roots of SAA can be traced back to M-estimation~\cite{1964hub,2014rus,2002steboo} from classical statistics.) SAA, unlike SG, is not an algorithm per se but a ``framework'' --- in SAA, a \emph{sample-path approximation} $f_m$ is ``constructed'' using an implicitly generated fixed sample $Y_1, Y_2, \ldots, Y_m$ in an attempt to approximate the objective $f$, and on which a chosen deterministic optimization algorithm is executed to high accuracy to obtain an estimator $X_m^*$ of a solution to~\eqref{master}. Traditionally, investigations in SAA focus on inferential questions about relevant SAA estimators, e.g., $f_m(X_m^*)$, $X_m^*$, and $f(X_m^*)$, under different structural conditions of the objective $f$ and its gradient $\nabla f$, and the nature of the stochastic oracle at hand. 

Unsurprisingly,  SAA solution estimators tend to be ``biased towards'' the chosen sample~\cite{1999makmorwoo}, tend to have high per-iteration computational cost due primarily to using a high sample size for approximating $f$, but also tend to enjoy superior iteration complexity (simply because the sampling effort as measured by the number of calls to the oracle are not recorded in the iteration complexity calculation). A crucial practical advantage that SAA enjoys over SG is that well-established deterministic optimization algorithms such as L-BFGS~\cite{nocedalbook} and Newton-CG~\cite{nocedalbook} can simply be ``plugged in'' as a solver in the service of optimizing $f_m$. 

\subsection{Retrospective Approximation (RA)}
The framework we present here, called ``Retrospective Approximation (RA)'' is designed to combine advantageous elements within SG and SAA, much like some ideas that have appeared in the literature over the last two decades~\cite{2016jalsha,2003tit,2021jalnedshayou,2009denfer,2015shabla}. RA and similar ideas  have also appeared before, in the context of stochastic programming problems~\cite{2021passon,2008polroy}, and stochastic root-finding problems~\cite{2001chesch,2010pas}. RA is organized into outer and inner iterations, where during each outer iteration, the incumbent iterate is updated by executing a deterministic optimization algorithm on a sample-path approximation of the objective, for a specified number of inner iterations. 

Three aspects of RA are salient. First, unlike in SA where a \textit{single} update of the iterate is performed using each sample generated for function and gradient approximation, RA allows \textit{multiple} updates of the iterate before considering a new sample for function and gradient approximation. Through such design, and assuming the use of common random numbers~(see Section~3), RA enjoys SAA's facility to ``plug-in'' existing quasi-Newton deterministic optimization algorithms such as L-BFGS and Newton-CG for performing the iterate updates. Such ability to directly use existing solvers avoids complex redesign of optimization algorithmic components such as steplength selection and Hessian estimation, and the necessity of tuning hyper-parameters, a procedure that is often necessary in SG~\cite{2017cursch}. 

Second, RA deviates from traditional SAA in that it explicitly considers the question of how many updates to perform before considering a new sample to approximate the function and its gradient. Comparing the statistical benefit derived from one additional iterate update versus a sample update is essentially a formal treatment of the \emph{estimation error versus optimization error} question discussed in~\cite{bottou2018optimization}, but within a sequential adaptive SAA context. RA answers this question by incorporating an adaptive stopping tolerance which dictates when it might be advantageous to stop iterate updates in favor of a sample update. 

A third aspect that works in favor of RA is the potential reduction in the number of times a new approximation of $f$ (or its gradient) is constructed, savings that might be significant in machine learning settings where such approximations might incur nontrivial computational cost.

\begin{remark}
The recent work in~\cite{2021jalnedshayou} requires special mention due to its similar point of view with RA in attempting to incorporate a quasi-Newton scheme alongside increasing sample sizes in the service of solving a stochastic optimization algorithm. However, RA differs in two respects that we think are crucial for both theoretical and implementation standpoints: (i) as we clarify in Section~\ref{sec:ra}, RA \emph{fixes} the sample between outer iterations, allowing a deterministic quasi-Newton solver to retain and exploit structure in the resulting deterministic optimization problem. In contrast, the quasi-Newton scheme in~\cite{2021jalnedshayou} relies on stochastic updating of the Hessian matrix obtained using \emph{differing} samples over iterations; and (ii) the deterministic solver being used in the outer iteration of RA is stopped using an adaptive error tolerance rule that judges whether it might be statistically efficient to continue optimizing the current sample-path problem.
\end{remark}

\subsection{Summary of Contributions}
RA's thrust is leveraging existing \emph{``off-the-shelf'' deterministic} optimization algorithms toward an implementable framework that reduces the need for hyper-parameter tuning. Our aim is not identifying the weakest assumptions in comparison with other popular frameworks such as SG, but instead to illustrate and clarify the principles that underlie RA. (For detailed discussion on assumptions in the context of SG, see~\cite{2020kharic,2021derkas,2022patzhatia}.) The key aspects of the paper are summarized next.
\begin{enumerate}
    \item[(a)] We propose Retrospective Approximation (RA) as a framework to solve smooth nonconvex stochastic optimization problems using a stochastic oracle. RA uses a chosen 
    solver, e.g., L-BFGS or Newton-CG, to solve a sequence of sampled deterministic optimization problems to a specified adaptive tolerance chosen to trade-off estimation error and optimization error. RA can also be viewed differently, as a generalized form of SG that incorporates the potential for multiple solution updates between sample updates.
    \item[(b)] We demonstrate that almost sure and $L_1$ consistency of iterates generated by RA are preserved under especially weak conditions. We characterize the iteration and oracle complexity when linear and sub-linear solvers are used within RA, and identify a practical solver termination criterion leading to optimal complexity rates. 
    \item[(c)] Our general premise is that controlling the \emph{relative error} of the observed gradients, as opposed to their \emph{absolute error}, is sufficient to ensure desirable behavior of algorithms for stochastic optimization. This is also necessary in many application domains as the underlying stochasticity very often comes in the form of \emph{multiplicative} noise or they depend on the gradient of the function. Accordingly, our entire mathematical treatment  assumes a uniform law on the relative error of the observed gradients in contrast with the more typical assumption of a uniform law on the absolute error~\cite{2006barjormca,2011botbos,2018meibaimon,2019mokozdjad,2009shadenrus,1995vap}. As we demonstrate in Section~\ref{sec:examples}, such treatment results in increased generality, allowing to subsume important problem classes otherwise excluded by a uniform law on the absolute error.
    \item[(d)] The objective $f$ is smooth and nonconvex and to reflect the effect of curvature at all first-order critical points of $f$, we establish a weak convergence result (a certain ``random CLT'') on the sequence of RA iterates, and the sequence of true gradients observed at RA's iterates. \item[(e)] Extensive numerical testing suggests that RA, due to its facility of incorporating established deterministic solvers without modification, seems to be able to dispense with hyper-parameter tuning. In particular, our numerical results demonstrate that higher-order methods when implemented within the RA framework are competitive to well-tuned SG methods. (All numerical tests in Section~\ref{sec:numerical} were performed using RA embedded with an L-BFGS line search solver.)  
    \item[(f)] Our use of relative error as the principal mathematical object combined with a reliance on some powerful recent results from empirical process theory~\cite{2006gee, 2016ginnic,1996vanwel} might present a standard route to generalizing various results on empirical risk minimization that have appeared in the machine learning literature.   
\end{enumerate}

\subsection{Paper Organization}

The paper is organized into nine sections plus an appendix. 
Section $2$ provides background, notation, and common assumptions used in the paper and Section $3$ describes the RA framework. Sections~4--8 develop the main ideas of the paper: Section~4 presents almost sure and $L_1$ consistency results for RA, Section~5 presents iteration and oracle complexity results, Section~6 presents 
four problem classes subsumed by RA, and  Section $7$ provides the Central Limit Theorems (CLT) for our RA. Section~8 presents some numerical experiments illustrating our theory. Section~9 gives some concluding remarks. The appendix provides 
the proof and an example for the complexity result for the weighted relative error class and also some essential statements from empirical process theory relevant to our RA setting.

\section{MATHEMATICAL PRELIMINARIES}\label{sec:notation} In what follows, we list notation and define key ideas that are repeatedly invoked in the paper.
 
\subsection{Important Notations} (i) If $x \in \mathbb{R}^d$, then $x_i, i = 1,2,\ldots,d$ refers to the $i$-th coordinate of $x$, implying that $x= (x_1,x_2,\ldots,x_d)$. (ii) For $x,y \in \mathbb{R}$, $a \vee b := \max(a,b)$ refers to the maximum of $a$ and $b$, and $a \wedge b := \min(a,b)$ refers to the minimum of $a$ and $b$. (iii) $\|x\|_p := \left( \sum_{i=1}^d |x_i|^p\right)^{\frac{1}{p}}$ refers to the $L_p$ norm of $x \in \mathbb{R}^d$. As usual, we use the special notation $\|x\|$ for the $p=2$ case to refer to the $L_2$ norm of $x \in \mathbb{R}^d$. (iv) For a random sequence $\{X_n, n \geq 1\}$, we write $X_n \as X$ to refer to almost sure convergence (also known as convergence with probability one), $X_n \inP X$ to refer to convergence in probability, and $X_n \inD X$ to refer to convergence in distribution (also known as weak convergence).  (v) For positive valued sequences $\{a_n, n \geq 1\}$ and $\{b_n, n \geq 1\}$, we say $a_n \sim b_n$ to mean $a_n/b_n \to 1$ as $n \to \infty$. (vi) The function $\delta_{ij} = 1$ if $i=j$ and $0$ otherwise refers to the Kronecker delta function.

\subsection{Important Definitions}

\begin{definition}[first-order critical points]\label{defn:criticalpts} The set $\mathcal{X}^*$ refers to the set of first-order critical points of the objective $f$, that is, 
\begin{equation}\label{critical.points}
\mathcal{X}^* := \{x \in \mathbb{R}^d: \|\nabla f(x) \| =0\}.
\end{equation}
\end{definition}
\begin{definition}[unbiased first-order stochastic oracle] We follow the oracle model detailed in~\cite{2012agaetal,1983nemyud} within a sample-path context. So, by an unbiased stochastic first-order oracle associated with an objective function $f$, we are referring to a random map $$(x,f,Y) \mapsto  (F(x,Y), \nabla F(x,Y)),$$ where $\mathbb{E}[F(x,Y)] = f(x)$ and $\mathbb{E}[\nabla F(x,Y)] = \nabla f(x)$. An analogous definition holds when considering an unbiased first-order stochastic oracle within a filtered probability space. 
\end{definition}
\begin{definition}[globally convergent] An iterative algorithm for identifying a first-order critical point of the function $f: \mathbb{R}^d \to \mathbb{R}$ is said to be \emph{globally convergent} if, starting with any initial guess $x_0$, the sequence of iterates $\{x_k, n \geq 1\}$ generated by the algorithm has gradient norms converging to zero, that is, \begin{equation}\label{globalconv}\|\nabla f(x_k)\| \to 0 \mbox{ as } k \to \infty.\end{equation} The iterative algorithm is \emph{locally convergent} if~\eqref{globalconv} holds only when the initial guess $x_0$ is such that $\| \nabla f(x_0)\|$ is small enough.
\end{definition}

\begin{definition}[linearly and sub-linearly convergent solver]\label{defn:solverspeed} Suppose an iterative solver generates iterates $\{x_k, k \geq 1\}$ starting at the initial guess $x_0$ when identifying a first-order critical point of a function $f: \mathbb{R}^d \to \mathbb{R}$. Such a solver is said to exhibit \emph{$(\rho,c)$-linear convergence} ($\rho < 1, c < \infty$) if the gradient norm at the iterate $x_k$ obtained after executing $n$ steps of the algorithm satisfies \begin{equation}\label{linconv}\| \nabla f(x_k) \|  \leq c \, \| \nabla f(x_0)\|\, \rho^k. \end{equation}  
A solver is said to exhibit \emph{$(\frac{1}{2},c)$-sublinear convergence} 
(with $c < \infty$) when executed on a function $f$ bounded from below, that is, $f^* = \inf\{x \in \mathbb{R}^d: f(x)\} > -\infty$, the gradient norm at the iterate $x_k$ obtained after executing $k$ steps of the algorithm satisfies \begin{equation}\label{sublinconv}\| \nabla f(x_k) \|  \leq c \, (f(x_0) - f(x_k))^{1/2} \, k^{-1/2} \leq c \, (f(x_0) - f^*)^{1/2} \, k^{-1/2}. \end{equation} 
\end{definition}

\begin{remark}\label{rem:sublindef}
The definition in~\eqref{sublinconv} expresses the progress made by a solver (in identifying a first-order critical point) in terms of the number of steps $k$ and the initial optimality gap $f(x_0)-f^*$. To further understand the structure in~\eqref{sublinconv}, consider the gradient descent recursion $$x_{t+1} = x_t - \eta \nabla f(x_t), t = 0,1,\ldots; \quad x_k^* := \arg\min\{\|\nabla f(x_t) \|: 1 \leq t \leq k\} $$ when executed with step size $\eta \leq L^{-1}$ on an $L$-smooth, lower bounded function $f:\mathbb{R}^d \to \mathbb{R}$, $- \infty < f^* := \inf_{x \in \mathbb{R}^d} f(x)$. It is well-known~\cite[pp. 28]{2004nes} that the sequence $\{x_k^*, k \geq 1\}$ exhibits the global sub-linear convergence rate \begin{equation}\label{graddescentconvrate} \|\nabla f(x^*_k)\| \leq  \left(\frac{2}{\eta}\right)^{1/2} \left( f(x_0) - f^*\right)^{1/2}k^{-1/2},
\end{equation} thereby satisfying~\eqref{sublinconv} with $c=\sqrt{\frac{2}{\eta}}$. We present oracle complexity results based on~\eqref{linconv} and~\eqref{sublinconv}.
\end{remark}

\begin{definition}[Empirical Measure and Empirical Process]\label{defn:emp} 
Let $\{Y_j, 1 \leq j \leq n\}$ be a sequence of independent and identically distributed (IID) copies of a random variable $Y$ having distribution $P$ on some appropriate measurable space $(\mathcal{Y}, \mathcal{A})$. The \emph{empirical measure} $P_n$ associated with $\{Y_j, j \geq 1\}$ is 
$$
P_n(A) := n^{-1} \sum_{j=1}^n \mathbb{I}_A (Y_j),
$$
where $A \in \mathcal{A}$. 
Let $P(A) := P\left\{Y\in A\right\}$ which can also be written as $\mathbb{E}\mathbb{I}_A(Y)$.
We define the \emph{empirical process} as the following ``scaled'' version of the above:
$$
G_n(A) := \sqrt{n}\left(P_n(A) - P(A)\right).
$$
The above can also be defined for a general function $F$ of $Y$:
$$
P_nF := n^{-1}\sum_{j=1}^n F(Y_j),\quad
G_n(F) := \sqrt{n}(P_nF - PF),
\quad
\text{where $PF = \mathbb{E}F(Y)$}.
$$
\end{definition}

By the classical law of large numbers and central limit theorems, for any fixed $A\in \mathcal{A}$ and any fixed function $F$ of $Y$, we have $P_n(A)$ and $P_nF$ converging almost surely to $P(A)$ and $PF$, and 
$G_n(A)$ and $G_nF$ converging in distribution to normal random variables with mean $0$ and variances $P(A)-P(A)^2$ and $\mathbb{E}(F-PF)^2$.
However, for many applications, we would like the above limit theorems to hold true \emph{uniformly} over some class of subsets $A$ and functions $F$.

The above concepts are very relevant to our current paper as our goal is exactly to minimize the empirical version of the random function $F$, with the ability to observe
$$
\nabla f_n(x) := n^{-1} \sum_{j=1}^n \nabla F(x,Y_j)
\,\,\,\left(=P_n\nabla F(x,\cdot)\right).
$$ 
Note that now, instead of indexed by a set $A$, we are dealing with the class of functions indexed by $x\in\mathcal{R}^n$, $\{\nabla F(x,\cdot): x \in \mathbb{R}^d\}$. Analogously, the ``scaled'' stochastic process 
$$
\Big\{\sqrt{n}\big(\nabla f_n(x) - \nabla f(x)\big), x  \in \mathbb{R}^d\Big\}
\quad\text{where $f(x) = \mathbb{E}F(x,Y)$}
$$ is called the \emph{empirical process}.

\section{RA DESCRIPTION}\label{sec:ra}

In this section, we define RA more precisely and as it is applied to the type of stochastic optimization problems we consider in this paper. 

Let $\{\tilde{Y}_k, k \geq 1\}$ be an IID ($\mathcal{Y}$-valued) random sequence on $(\Omega, \mathcal{F},P)$. To precisely define a sequence of empirical mean function and gradient estimators, let's construct an increasing sequence of sigma algebras $\mathcal{F}_1 \subseteq \mathcal{F}_2 \subseteq \cdots $ as follows. Let $M_1$ be a constant, define 
\begin{align}\label{sigmaalg} 
\mathcal{F}_1 &:= \sigma(\tilde{Y}_j \vert j \leq M_1); && \text{and}\,\,\,M_2 \in \mathcal{F}_1; \nonumber \\ \mathcal{F}_2 &:= \sigma(\tilde{Y}_j \vert j \leq M_1 + M_2); && \text{and}\,\,\,M_3 \in \mathcal{F}_2; \nonumber \\ \vdots & && \nonumber \\  
\mathcal{F}_k &:= \sigma(\tilde{Y}_j \vert j \leq \sum_{j=1}^{k} M_j); && \text{and}\,\,\,M_{k+1} \in \mathcal{F}_k.  \end{align} Let $\tilde{Y}$ be a copy identically distributed with $\{\tilde{Y}_k, k \geq 1\}$.  Define the \emph{sample-path function} $f_m$ and the \emph{sample-path gradient} having sample size $m$ as follows:

\begin{align}\label{samplepathfnandgrad}
 f_{m}(x) := \frac{1}{m}\sum_{j=1}^{m} F(x,\tilde{Y}_j); \quad \nabla f_m(x) := \frac{1}{m}\sum_{j=1}^m \nabla F(x,\tilde{Y}_j) \quad x \in \mathbb{R}^d.  
 \end{align} 
 Thus we see that sample-path function $f_{m}(\cdot)$ and the sample-path gradient $\nabla f_m$ appearing in~\eqref{samplepathfnandgrad} are constructed as the sample mean of the random functions  $F(\cdot,\tilde{Y}_j), j =1,2,\ldots,m$ and their gradients $\nabla F(\cdot,\tilde{Y}_j), j =1,2,\ldots,m$, respectively, where $\tilde{Y}_j,j=1,2,\ldots,m$ are identically distributed copies of a random object $\tilde{Y}$. The random object $\tilde{Y}$ will never enter the stage of our presentation directly, but we assume that the gradients $\nabla F(\cdot,\tilde{Y})$ are unbiased with respect to $\nabla f(\cdot)$, that is, $$\mathbb{E}\left[ \nabla F(x,\tilde{Y}) \right] = \nabla f(x), \quad x \in \mathbb{R}^d.$$ Due to the above construction of the filtration $(\Omega,\mathcal{F},(\mathcal{F}_k)_{k \geq 1},P)$, notions of a sample-path function and a sample-path gradient implicitly constructed with an adapted sample size $M_k \in \mathcal{F}_{k-1}$ are now well defined:
 \begin{align}\label{relabel}
 f_{M_k}(x) := \frac{1}{M_k}\sum_{j=1}^{M_k} F(x,Y_j); \quad \nabla f_{M_k}(x) := \frac{1}{M_k}\sum_{j=1}^{M_k} \nabla F(x,Y_j) \quad x \in \mathbb{R}^d,
 \end{align} where in~\eqref{relabel}, we have abused notation to obtain the convenient relabeling $Y_j := \tilde{Y}_{\left(\sum_{i=1}^{k-1} M_i\right) + j}$. This definition allows for the sample sizes $M_k$ to be stochastic, unlike in standard minibatch stochastic gradient (SG) methods \cite{bottou2018optimization}. While requiring elaborate filtration machinery for precise definition,~\eqref{relabel} conveys the simple and intuitive idea that the function and gradient estimate at any point $x$ during iteration $k$ are each obtained as the sample mean of $M_k$ unbiased observations. Furthermore, while we develop much of the theory that follows assuming that the sequence $\{\tilde{Y}_k, k \geq 1\}$ is IID, a relaxation of the independence assumption is possible using the now standard machinery from empirical process theory~\cite{1996vanwel, 2016ginnic, 2006gee}. Such relaxation will subsume important contexts where $\{\tilde{Y}_k, k \geq 1\}$ constitutes observations from a real system, e.g., a queuing network operating at equilibrium.  
 
 We assume throughout that we have at our disposal a \emph{globally convergent} (see Definition~\ref{globalconv}) ``iterative solver'' that we generically call Method $\mathfrak{M}$, capable of solving the sample-path optimization problem: \begin{align*}\label{samplepathproblem}
 \mbox{minimize: } & f_{M_k}(x) \nonumber \\ \mbox{ subject to: } & x \in \mathbb{R}^d. \tag{$Q_{M_k}$}\end{align*} More precisely, we assume having at our disposal Method~$\mathfrak{M}$ which, given any $\epsilon >0$ and starting with any initial guess, identifies a point $x \in \mathbb{R}^d$ satisfying $\| \nabla f_{M_k}(x)\| \leq \epsilon.$ Two examples of the Method $\mathfrak{M}$ that are of particular interest are line search and trust-region procedures~\cite{nocedalbook} executed with L-BFGS~\cite{nocedalbook} for Hessian estimation.
 
 Next we give the listing of the RA Algorithm.
 
 \begin{algorithm}[H]
\caption{RA} 
\label{alg:RA}
\begin{algorithmic}[1]
\STATE{{\bf Input:} (i) initial guess $X_0$; (ii) procedure to update $\{M_k, k \geq 1\}$; (iii) procedure to update $\{\epsilon_k, k \geq 1\}$; (iv) Method $\mathfrak{M}$, e.g., line search or trust-region with L-BFGS.} 
\vspace{2mm}
 \FOR{$k = 1,2,\cdots$}
 \vspace{2mm}
 \STATE{{Set sample size and error tolerance:} Choose $M_k \in \mathcal{F}_{k-1}$ and $\epsilon_k \in \mathcal{F}_{k-1}$.} 
        \STATE{{Execute Method~$\mathfrak{M}$ for $N_k$ steps:} With Method $\mathfrak{M}$ and $X_{k-1}$ as ``warm-start,'' execute $N_k$ steps to obtain `inner iterates'' $X_{k-1,t}, t=0,1,2,\ldots,N_k$ satisfying $$\|\nabla f_{M_k}(X_{k-1,N_k})\| \leq \epsilon_k, \quad \mbox{ e.g., } \epsilon_k := \Gamma_k \|\nabla f_{M_k}(X_{k-1,0})\|, \mbox{  where } \Gamma_k > 0$$}
\STATE{{Set the $k$-th outer iterate:} $X_k = X_{k-1,N_k}$.} 
\ENDFOR
\end{algorithmic}
\end{algorithm}

As seen in the algorithm listing, Algorithm~RA is an iterative framework that is organized into outer and inner iterations. During the $k$-th (outer) iteration, Method~$\mathfrak{M}$ uses the $(k-1)$-th outer iterate as the initial guess to execute $N_k$ steps and solve the sample-path problem~\eqref{samplepathproblem} to within adapted tolerance $\epsilon_k \in \mathcal{F}_{k-1}$. Let $\{X_k, k \geq 1\}$ denote the (outer) iterates, and $X_{k-1,t}, t = 0,1,2,\ldots,N_k$ the inner iterates generated by Method~$\mathfrak{M}$ during iteration $k$. Since Method~$\mathfrak{M}$ uses $X_{k-1}$ as the initial guess during the $k$-th iteration and terminates after $N_k$ steps when the norm of the sample-path gradient drops below $\epsilon_k$, we see that the following hold: $$X_{k-1,0} := X_{k-1}; \quad X_k := X_{k-1,N_k}; \quad \| \nabla f_{M_k}(X_k) \| \leq \epsilon_k.$$   

We now briefly discuss four aspects of RA that are of vital theoretical and practical importance.
\subsection{Choice of Method~$\mathfrak{M}$.} One of the key advantages of RA is that its design naturally allows ``directly plugging-in'' a highly successful deterministic optimization solver (generically called Method~$\mathfrak{M}$ in Algorithm RA) in the service of solving a stochastic optimization problem. While Method~$\mathfrak{M}$ can be any globally convergent procedure for identifying a first-order critical point in~\eqref{samplepathproblem}, our particular interest is a quasi-Newton method such as line search or trust-region with L-BFGS~\cite{nocedalbook} for Hessian estimation where the updates to BFGS matrices are skipped if the \emph{curvature condition} is not satisfied (see \cite{BollapragadaICML18}). 

Step~4 in Algorithm~RA assumes that the Method~$\mathfrak{M}$ is able to solve~\eqref{samplepathproblem} to any prescribed level of tolerance $\epsilon_k$, that is, Method~$\mathfrak{M}$ is able to identify a point $X_k$ such that $\|\nabla f_{M_k}(X_k)\| \leq \epsilon_k$. Such assumption is valid for standard quasi-Newton methods such as line search or trust-region~\cite{nocedalbook} with L-BFGS~\cite{nocedalbook} if the random functions $F(\cdot,Y)$ that comprise the sample-path function are $L(Y)$-smooth, that is, $$
\|\nabla F(x,Y) - \nabla F(y,Y)\| \leq L(Y)\|x-y\|, \quad \forall x,y  \in \mathbb{R}^d,
\quad\text{with $\mathbb{E}[L^2(Y)] < \infty.$}
$$
It is well-known~\cite{2015kimpashen,2009shadenrus} that when $F(\cdot,Y)$ is $L(Y)$-smooth and $\mathbb{E}[L^2(Y)] < \infty$, smoothness is transferred to the objective $f$, that is, the function $f(x) = \mathbb{E}[F(x,Y)], x \in \mathbb{R}^d$ is $\mathbb{E}[L]$-smooth.

\subsection{Choice of Sample Size, Error Tolerance, Termination Criterion.} Recall that RA executes Method~$\mathfrak{M}$ on the sample-path problem~\eqref{samplepathproblem} during the $k$-th iteration. Since the sample-path problem ~\eqref{samplepathproblem} is only a proxy for the true problem~\eqref{master}, it stands to reason that Method~$\mathfrak{M}$ not be executed ad infinitum but instead terminated if further solving is deemed to produce no statistical gains. RA formalizes this idea by introducing the notion of an error tolerance $\epsilon_k$ in Step 4 of RA: during the $k$-th iteration Method~$\mathfrak{M}$ executes until it finds an iterate $X_{k-1,N_k}$ satisfying $\|\nabla f_{M_k}(X_{k-1,N_k}\| \leq \epsilon_k$.

The choice of the sample size $M_k$ and the tolerance $\epsilon_k$ are central questions in our later analysis, but a common and convenient choice that is consistent with optimization practice is to choose $$\epsilon_k := \Gamma_k \|\nabla f_{M_k}(X_{k-1,0})\| = \Gamma_k \|\nabla f_{M_k}(X_{k-1})\| \mbox{ and } M_k = c_1M_{k-1}, \mbox{ e.g. } c_1=1.05$$ implying that the Method~$\mathfrak{M}$ executes as long as the norm of the sample gradient at the incumbent iterate remains above that at the initial guess by at least $\Gamma_k$, and the sample sizes across outer iterations $k$ are increased by a fixed percentage $100c_1$. We calculate the oracle complexity resulting from this choice in Remark~\ref{rm:workcomp_alt}.

\subsection{Warm Starts.} RA employs ``warm starts,'' that is, when executing Method~$\mathfrak{M}$, the solution to the previous iteration becomes the initial guess for the subsequent iteration. The use of warm starts is more than just an implementation heuristic and features explicitly in the oracle complexity calculations of RA. Specifically, warm starts reduce the number of inner iterations performed, and such reduction is featured prominently when the incumbent sample size is high. 

\subsection{Common Random Numbers} The effectiveness of RA, particularly its use of a quasi-Newton solver to solve sample-path problems, is crucially dependent on the use of what is called \emph{common random numbers} (CRN)~\cite{nel2013}. Formally, notice that RA uses $Y_1,Y_2, \ldots, Y_{M_k}$ to generate the sample-path objective function $f_{M_k}$, that is, the ``sample remains fixed'' as the search procedure Method~$\mathfrak{M}$ in Step 4 of RA is executed on $f_{M_k}(\cdot)$ for $N_k$ steps towards solving the sample-path problem~\eqref{samplepathproblem}. CRN should be seen as a mechanism that preserves any sample-path structure that may be present, e.g., strong convexity, in the process allowing faster convergence of the deterministic solver.

\section{CONSISTENCY}

In this section, we present results that characterize the convergence behavior of the sequence $\{\| \nabla f(X_k) \|, k \geq 1\}$ of true gradient norms evaluated at the iterates generated by RA. To preserve generality and applicability, we will advocate that the relative error of the gradient estimator is a principal object dictating algorithmic behavior. More precisely, the \emph{relative error class} $\mathcal{R}$ is defined as: \begin{equation}\label{relerror} \mathcal{R} : = \{R(x,Y), x \in \mathbb{R}^d, Y \in \mathcal{Y}\}; \quad R(x,Y) := \frac{\nabla F(x,Y) - \nabla f(x)}{\| \nabla f(x) \| + c_0}; \quad c_0 >0.
\end{equation} 
The constant $c_0$ is to avoid the nuisance associated with a vanishing gradient in the stationary set $\mathcal{X}^*$ (defined through \eqref{critical.points}). 
The class $\mathcal{R}$ should be interpreted as the collection of the relative error functions $R(x,\cdot)$, labeled by $x \in \mathbb{R}^d$.
The empirical version of the above is given by
\begin{equation}\label{empavg} 
R_m(x) := \frac{ \nabla f_m(x) - \nabla f(x) }{\| \nabla f(x) \| + c_0}, \quad x \in \mathbb{R}^d,
\quad\text{where}\quad
\nabla f_m(x) := m^{-1} \sum_{j=1}^m \nabla F(x,Y_j).
\end{equation}
In addition, we define the \emph{envelope} and the empirical envelope as
\begin{equation}\label{defn:env} G = G(Y) := \sup \left\{\|R(x,Y)\|: x \in \mathbb{R}^d\right\},\end{equation}  and
\begin{equation}\label{eq:Gm}
G_m(Y) := \sup\{\|R_m(x)\|: x \in \mathbb{R}^d\}.
\end{equation}


The example problem classes of Section~\ref{sec:examples} serve to support our claim that the \emph{relative error class} affords great generality, and that \emph{relative error} but not \emph{absolute error} of the gradient estimator is a fundamental object that decides consistency and convergence rates in stochastic optimization. 
This is a reflection of the fact that in many applications, the noise is multiplicative, very often dependent on the gradient $\nabla f(X)$ which can be unbounded in the overall simulation space. Hence it is more important to consider $\displaystyle \frac{\|\text{noise}\|}{\|\nabla f(X)\|}$ rather than simply $\|\text{noise}\|$.
The simple but widely used \emph{multiplicative error} classes
\eqref{mulerrclass0} (in Theorem~\ref{thm:itercomp}), \eqref{mulerrclass} (in Section~\ref{ex:multerrclass}), and \eqref{mulerrclass2} (in Section ~\ref{PECEx}) all illustrate this point well. These examples satisfy  stipulations on the uniform law of the relative error for strong consistency but violates the corresponding routinely assumed requirements for consistency on the absolute error. 

\subsection{Main Results on Consistency}
In what follows we demonstrate RA's strong consistency property by combining conditions on the behavior of the termination criterion sequence $\{\epsilon_k, k \geq 1\}$ and the sample size sequence $\{M_k, k \geq 1\}$ with a uniform law on the relative error class. 

\begin{theorem}[Strong Consistency of RA]\label{thm:consistency} Suppose that the relative error class $\mathcal{R}$ defined in~\eqref{relerror} satisfies the uniform law \begin{equation}\label{ass:RULLN} \sup_{x \in \mathbb{R}^d} \|R_{m}(x)\| \as 0 \mbox{ as } m \to \infty, \tag{R-ULLN}\end{equation} and that \begin{equation} \label{cond:errtolzero} \epsilon_k \,\, \as \,\,  0 \mbox{ and } M_k \,\, \as \,\,  \infty \mbox{ as } k \to \infty. \end{equation}  Then \begin{equation} \label{RAwp1conv}  \| \nabla f (X_k) \| \,\, \as \,\, 0 \mbox{ as } k \to \infty.\end{equation}
\end{theorem}

Theorem~\ref{thm:consistency} gives sufficient conditions for consistency in probability and  with probability one under conditions on the relative error class $\mathcal{R}$ that match standard assumptions in stochastic optimization~\cite[Assumption~4.3]{bottou2018optimization} and in the language of empirical process theory~\cite{2006gee} (see Section~\ref{sec:examples} for example problem classes). 
A third type of consistency that is more routinely employed in the optimization literature is what we call $L_1$ consistency, that is, $\mathbb{E}\left[\|\nabla f (X_k) \| \right] \to 0.$ It is important that none of these modes of convergence provide guarantees on a solution estimator $X_k$ obtained through a suitable stopping heuristic. While important, this question of assessing the quality of a solution obtained using sequential stopping lies outside the scope of this paper --- see~\cite{2022pat} for more.  

As will become evident shortly through the proof of Theorem~\ref{thm:l1cons}, in addition to the uniform convergence \eqref{ass:RULLN} of the relative error class $\mathcal{R}$, $L_1$ consistency follows if the iterates are such that $\mathbb{E}\left[\|\nabla f (X_k) \|^2 \, \vert \, \mathcal{F}_{k-1}\right]$ remains bounded, implying in effect that the variance in the iterates is not so high that the norm of the true gradient evaluated at the iterates does not diverge in expectation. As the examples in Section~\ref{sec:examples} will illustrate, our use of relative error demonstrates that this condition is widely satisfied. Similar to Theorem~\ref{thm:consistency}(a) for convergence in probability, Theorem~\ref{thm:l1cons} also imposes the condition $\epsilon_k \to 0$ on the error tolerance sequence $\{\epsilon_k, k \geq 1\}$; this is in contrast to Theorem~\ref{thm:consistency}(b) for almost sure convergence,  which imposes the slightly more stringent condition $\sum_k M_k^{-1} < \infty$ on the growth of the sample size sequence $\{M_k, k \geq 1\}$.

\begin{theorem}[$L_1$ Consistency of RA]\label{thm:l1cons}  Suppose the sample size sequence $\{M_k, k \geq 1\}$ and the error-tolerance sequence $\{\epsilon_k, k \geq 1\}$ are chosen so that $$M_k \as \infty; \quad \mathbb{E}\left[\epsilon_k \, \vert \, \mathcal{F}_{k-1}\right] \as 0,$$ and the relative error class $\mathcal{R}$ defined in~\eqref{relerror} satisfies the uniform law in~\eqref{ass:RULLN} along with $\mathbb{E}[G^2(Y)] < \infty$. Suppose also that the iterates $\{X_k, k \geq 1\}$ are such that \begin{equation}\label{graditerbd} \sup_k \, \mathbb{E}[\| \nabla f(X_k) \|^2 \, \vert \, \mathcal{F}_{k-1}] =: \ell^2_{2,\infty} < \infty. \tag{\mbox{Bd-Grad}} \end{equation} Then as $k \to \infty$, \begin{equation}\label{l2bdgradfconv}\mathbb{E}\left[ \| \nabla f(X_k) \| \, \vert \, \mathcal{F}_{k-1} \right] \as 0\end{equation} 
and also
\begin{equation}\label{l1conv}
  \mathbb{E}\big[ \| \nabla f(X_k) \| \, \big] \longrightarrow 0.\end{equation}\end{theorem}

\begin{remark}\label{rem:equivalence} Theorem~\ref{thm:consistency} prominently assumes a uniform law on the relative error class. As we shall see later, verifying this condition is fairly direct in the context of optimization, in particular for explicit examples. On the other hand, it is instructive to relate this condition to a more fundamental property about the underlying stochasticity of the oracle. Qualitatively, this relation touches upon the 
``regularity'' or ``irregularity'' of the stochastic functions
$F(x,Y)$ used by the oracle to approximate the function $f$ to be optimized. This has been studied extensively in the field of empirical processes. The quantitative measure is given by the entropy of the class of functions $\left\{F(x,Y), x\in\mathcal{R}^d\right\}$. Intuitively, if the entropy does not grow too fast as the fineness of approximation increases, the approximating functions $f_n$ can represent the true function $f$ more faithfully so that minimizing $f_n$ can give a more trustworthy answer. In fact, there is a more or less equivalence between the uniform law and the sub-linear growth of the entropy function.
We refer to \cite{1996vanwel,2006gee, 2016ginnic} for comprehensive theory. Some of the key results will be illustrated in the appendix.
\end{remark}

\subsection{Proofs of Theorem~\ref{thm:consistency} and Theorem~\ref{thm:l1cons}.}
\begin{proof}{Proof of Theorem~\ref{thm:consistency}.}  {Define $$\tilde{t} := \frac{t}{2c_0 + t} = 1 - \frac{2c_0}{2c_0 + t}, \quad t \in (0,1).$$ Since~\eqref{ass:RULLN} holds along with $M_k \as \infty$ and $\epsilon \as 0$, \begin{align}\label{pcons}  \lim_{K \to \infty} P\left( \bigcup_{k \geq K} \, \sup_{x \in \mathbb{R}^d} \| R_{M_k}(x) \| > \tilde{t} \right)  = 0 \mbox{ and } \lim_{K \to \infty} P\left( \bigcup_{k \geq K} \, \epsilon_k > c_0 \, \tilde{t} \right)  = 0.
\end{align}
Also, notice that there exists $K_1$ such that for $k \geq K_1$, \begin{align*}\label{tailsetequiv} \left\{ \| \nabla f(X_k) \| > t \right\} & \subseteq \left\{ (1 - \tilde{t})\| \nabla f(X_k)\| > \tilde{t}c_0 + \epsilon_k\right\} \nonumber \\   & \subseteq \left\{ (1 - \tilde{t})\| \nabla f(X_k)\| > \tilde{t}c_0 + \| \nabla f_{M_k}(X_k)\|\right\} \nonumber \\
& \subseteq \left\{ \| \nabla f_{M_k}(X_k) - \nabla f(X_k) \| > \tilde{t}\left(\| \nabla f(X_k) \| + c_0 \right) \right\} \nonumber \\
& \equiv \left\{ \| R_{M_k}(X_k) \| > \tilde{t} \right\} \cup \left\{ \epsilon_k > c_0\, \tilde{t} \right\},
\end{align*} and therefore, \begin{equation} \label{step1}\bigcup_{k \geq K_1} \left\{ \| \nabla f(X_k) \| > t \right\} \subseteq \bigcup_{k \geq K_1} \left\{ \| R_{M_k}(X_k) \| > \tilde{t} \right\} \cup \left\{ \epsilon_k > c_0\, \tilde{t} \right\}. \end{equation} Use~\eqref{pcons} and~\eqref{step1} to conclude that for any $\delta>0$, there exists $K_2=K_2(\delta)$ such that \begin{align}\label{step2} P\left(\bigcup_{k \geq K_2} \left\{ \| \nabla f(X_k) \| > t \right\}\right) &\leq P\left( \bigcup_{k \geq K_2} \| R_{M_k}(X_k) \| > \tilde{t} \right) + P\left( \bigcup_{k \geq K_2} \epsilon_k > c_0\, \tilde{t} \right) \nonumber \\ & \leq \frac{\delta}{2} + \frac{\delta}{2} = \delta, \end{align} implying that Theorem~\ref{thm:consistency} holds.}

\end{proof}

\begin{proof}{Proof of Theorem~\ref{thm:l1cons}.}
Recall the empirical-mean relative error $R_m(x)$ and the empirical-mean envelope as defined in~\eqref{empavg} and~\eqref{eq:Gm}. 
We see after some algebra that
\begin{align}
    \| \nabla f(X_k)\| &\leq \| \nabla f_{M_k}(X_k)\| + \| \nabla f_{M_k}(X_k) - \nabla f(X_k)\| \nonumber\\
   &\leq \| \nabla f_{M_k}(X_k)\| + R_{M_k}(X_k)(c_0 + \|\nabla f(X_k)\|) \nonumber \\
   &\leq \| \nabla f_{M_k}(X_k) \| + c_0 G_{M_k} + G_{M_k} \| \nabla f (X_k) \| \label{l1mainineq}
\end{align}
Recalling that $\|\nabla f_{M_k} (X_k)\| \leq \epsilon_k$, taking conditional expectations on~\eqref{l1mainineq}, and then applying Cauchy-Schwarz inequality, we get \begin{align} \label{expgenineqrelinit}\mathbb{E}\left[\| \nabla f(X_k) \| \, \vert \, \mathcal{F}_{k-1}\right] &\leq \mathbb{E}[\epsilon_k \, \vert \, \mathcal{F}_{k-1}] + c_0\mathbb{E}\left[G_{M_k} \, \vert \, \mathcal{F}_{k-1}\right] + \mathbb{E}\left[G_{M_k} \| \nabla f(X_k) \| \, \vert \, \mathcal{F}_{k-1}\right] \nonumber \\ &\leq \mathbb{E}[\epsilon_k \, \vert \, \mathcal{F}_{k-1}] + c_0\mathbb{E}\left[G_{M_k} \, \vert \, \mathcal{F}_{k-1}\right] + \left(\mathbb{E}\left[G^2_{M_k} \, \vert \, \mathcal{F}_{k-1}\right]\right)^{1/2} \left(\mathbb{E}\left[\|\nabla f(X_k)\|^2 \, \vert \, \mathcal{F}_{k-1}\right]\right)^{1/2} \nonumber \\ &\leq \mathbb{E}[\epsilon_k \, \vert \, \mathcal{F}_{k-1}] + c_0\mathbb{E}\left[G_{M_k} \, \vert \, \mathcal{F}_{k-1}\right] + \left(\mathbb{E}\left[G^2_{M_k} \, \vert \, \mathcal{F}_{k-1}\right]\right)^{1/2} \, \ell_{2,\infty},  \end{align} where the last inequality in~\eqref{expgenineqrelinit} follows from~\eqref{graditerbd}. However, since we have assumed that $\mathbb{E}[G^2(Y)]< \infty$ and that~\eqref{ass:RULLN} hold,  it follows that \begin{equation}\label{empenvl2zero}\mathbb{E}\left[G^2_{M_k} \, \vert \, \mathcal{F}_{k-1}\right] \as 0 \mbox{ and } \mathbb{E}\left[G_{M_k} \, \vert \, \mathcal{F}_{k-1}\right] \as 0.\end{equation}  Using~\eqref{empenvl2zero} and the assumption $\mathbb{E}\left[\epsilon_k \, \vert \, \mathcal{F}_{k-1}\right] \as 0$ in~\eqref{expgenineqrelinit}, we see that the assertion in~\eqref{l2bdgradfconv} holds. Statement \eqref{l1conv} follows from \eqref{graditerbd} and Lebesgue Dominated Convergence Theorem:
\[
\mathbb{E}[|\nabla f(X_k)|\,\vert\mathcal{F}_{k-1}]
\leq
\mathbb{E}[|\nabla f(X_k)|^2\,\vert\mathcal{F}_{k-1}] \leq \ell^2_{2,\infty} < \infty
\]
and hence 
$\mathbb{E}\big[\|\nabla f(X_k)\|\big]
=
\mathbb{E}\left[\mathbb{E}\big[\|\nabla f(X_k)\|\,\vert\mathcal{F}_{k-1}\big]\right]
\longrightarrow 0.
$
$\blacksquare$
\end{proof}

\section{COMPLEXITY CALCULATION}\label{sec:itercomp}

In this section, under a further strengthening of the uniform law assumption in~\eqref{ass:RULLN}, we provide results that quantify the non-asymptotic convergence rates (also known as \emph{complexity}) of RA when expressed in terms of the number of iterations (called \emph{iteration complexity}) and the expended number of oracle calls (called \emph{oracle complexity}). We treat each of these complexities in order. 

\subsection{Iteration Complexity}

It should be clear that both the error tolerance sequence $\{\epsilon_k, k \geq 1\}$ and the sample size sequence $\{M_k, k \geq 1\}$ interact to determine the iteration complexity. To understand these effects better, for all results that follow, we assume that for each $k \geq 1,$ \begin{equation}\label{errtolchoice} \epsilon_k := \frac{C_{2,k}}{M_k^{1/2}} \quad \mbox{ and } \quad  c_2 \geq C_{2,k} \in \mathcal{F}_{k-1}, \tag{C.2(a)} \end{equation} where $C_{2,k}$ is an $\mathcal{F}_{k-1}$-measurable random variable with the upper bound $c_2$. 
We are now ready to state the main iteration complexity result.

\begin{theorem}[Iteration Complexity]\label{thm:itercomp0} 
Let the postulates of Theorem~\ref{thm:l1cons} hold. Furthermore, let the empirical-mean envelope $G_m : = \sup_{x \in \mathbb{R}^d} \|R_{m}(x)\|$ exhibit ``CLT-scaling'', that is, there exists $\sigma < \infty$ such that,  \begin{equation}\label{emperrenv-clt0} \mathbb{E}\left[\left(\sqrt{M_k}\, G_{M_k}\right)^2 \, \vert \, \mathcal{F}_{k-1}\right] \leq \sigma^2; \quad M_k \in \mathcal{F}_{k-1}. \tag{CLT-sc}\end{equation} Then, for $k \geq 1$ we have that \begin{equation}\label{geniter0} \mathbb{E}\left[ \| \nabla f (X_k) \| \, \vert \, \mathcal{F}_{k-1} \right] \leq  \frac{c_2}{M_k^{1/2}} + \frac{(c_0\sigma + \ell_{2,\infty}\sigma)}{\sqrt{M_k}}. \end{equation} 

\noindent Choose the sample-size sequence $\{M_k, k \geq 1\}$ so that $M_1=m_1$ is some positive constant,   \begin{equation}\label{sampsizeseqchoice} M_k := C_{1,k} M_{k-1}; \quad  c_1 \leq C_{1,k} \in \mathcal{F}_{k-1}, \tag{C.1(a)}\end{equation} 
and the error tolerance as in~\eqref{errtolchoice}. Then for $k \geq 1,$ 
\begin{align} \label{itercomp}
    \mathbb{E} \left[\|\nabla f(X_k)\|\right] 
    \leq \bigg(\frac{1}{\sqrt{c_{1}}}\bigg)^{k} \bigg( \frac{(c_2 + c_0\sigma + \ell_{2,\infty}\sigma)}{\sqrt{m_1}} \bigg).
\end{align}
Furthermore, for all $k \geq K(\epsilon)$ where \begin{equation} \label{epsiloncomp0} K(\epsilon):= 1 + \frac{1}{\log \sqrt{c_1}}\left(\log \frac{1}{\epsilon} + \log \frac{c_2 + c_0\sigma + \ell_{2,\infty}\sigma}{\sqrt{m_1}}\right),\end{equation} we have that $$\mathbb{E}\left[ \| \nabla f (X_k) \| \, \vert \, \mathcal{F}_{k-1} \right] \leq \epsilon, \quad \epsilon > 0.$$
\end{theorem}  

\begin{proof}{Proof.} From~\eqref{expgenineqrelinit}, we have \begin{align} \label{expgenineqrel}\mathbb{E}\left[\| \nabla f(X_k) \| \, \vert \, \mathcal{F}_{k-1}\right] &\leq \mathbb{E}[\epsilon_k \, \vert \, \mathcal{F}_{k-1}] + c_0\mathbb{E}\left[G_{M_k} \, \vert \, \mathcal{F}_{k-1}\right] + \mathbb{E}\left[G_{M_k} \| \nabla f(X_k) \| \, \vert \, \mathcal{F}_{k-1}\right] \nonumber \\ &\leq \mathbb{E}[\epsilon_k \, \vert \, \mathcal{F}_{k-1}] + c_0\mathbb{E}\left[G_{M_k} \, \vert \, \mathcal{F}_{k-1}\right] + \left(\mathbb{E}\left[G^2_{M_k} \, \vert \, \mathcal{F}_{k-1}\right]\right)^{1/2} \left(\mathbb{E}\left[\|\nabla f(X_k)\|^2 \, \vert \, \mathcal{F}_{k-1}\right]\right)^{1/2}.  \end{align} From~\eqref{graditerbd}, we know that there exists $0 < \ell_{2,\infty} < \infty$ such that \begin{equation}\label{bditeratesagain} \mathbb{E}\left[\|\nabla f(X_k)\|^2 \, \vert \, \mathcal{F}_{k-1}\right] \leq \ell_{2,\infty}^2 < \infty. \end{equation} Use~\eqref{emperrenv-clt0},~\eqref{errtolchoice}, and~\eqref{bditeratesagain} to see that the assertion in~\eqref{geniter0} holds. To obtain~\eqref{itercomp}, recurse the upperbound in~\eqref{geniter0} backward. 
$\blacksquare$
\end{proof}

Several observations are salient. First, the iteration complexity in~\eqref{itercomp} expressed in terms of the ``outer iterations'' is linear, and looks like that in deterministic, smooth and strongly convex regimes. This is essentially because the work done by the solver in the inner iterations is not accounted. This is a key difference from classical iteration complexity rates using a method such as stochastic gradient descent. Next, we see that in addition to a uniform law on the relative error class that is needed for consistency of RA, iteration complexity entails a certain regularity of the envelope, expressed through the CLT-scaling assumption~\eqref{emperrenv-clt0}. This condition, loosely speaking, stipulates that the supremum norm of the relative error scales as $O(1/\sqrt{M_k})$ in expectation. 
Similar to the description in Remark~\ref{rem:equivalence}, such a CLT-scaling is related to some integrability condition of the entropy function. This condition can be expressed concisely in terms of the Koltchinskii-Pollard entropy function. See \cite[Theorems 3.5.1 and 3.5.4]{2016ginnic} and also the appendix for a brief exposition.
We expect the assumption in~\eqref{emperrenv-clt0} to hold widely, in particular for the examples in Section~\ref{sec:examples}
as it is made on the relative error envelope rather than the absolute error envelope.

Also, in deriving the iteration complexity, a more general relationship of the form $\epsilon_k = C_{2,k}/M_k^{p}, p>0$ could have been assumed in~\eqref{errtolchoice}. However, through calculations nearly identical to that in the proof of Theorem~\ref{thm:itercomp}, it can be shown in such a case that the fastest achievable iteration complexity by RA is $O(1/M_k^{p \wedge 1/2})$ suggesting the choice $p \geq 1/2$. We can further reason that one should choose $p=1/2$ since choosing $p>1/2$ amounts to ``oversolving," that is, leads to increased computational effort during each iteration $k$ without yielding a better complexity rate. We will make this idea of not ``oversolving" more precise in the ensuing section.

In the following, we give a complexity result that dispenses with the bounded iterates assumption in~\eqref{graditerbd} needed for $L_1$ consistency, while strengthening the assumption in~\eqref{ass:RULLN} using the notion of a \emph{weighted relative error class}. Combined with what have been presented, they illustrate a ``dichotomy'' between \emph{multiplicative} versus \emph{additive noise} in the oracle model. To illustrate this idea, suppose
\begin{equation}\label{mulerrclass0}
\nabla F(x,Y) = (I+{\cal A(Y)})\nabla f(x) + {\cal B(Y)}
\end{equation}
where ${\cal A(\cdot)}$ and $\cal B(\cdot)$ are some mean zero random variables. In this case, the relative error is given by
\[
R(x,Y) = \frac{{\cal A}(Y)\nabla f(x) + {\cal B}(Y)}
{||\nabla f(x)|| + c_0}
\]
Note that as long as $\cal A \not\equiv 0$, then $R(x,Y) \approx O(1)$ while if $\cal A \equiv 0$, we will have that 
$R(x,Y) = o(||\nabla f(x)||)$.
This motivates the following general statement which handles the latter but allows $\cal B$ to depend on $\nabla f(x)$ alas, to a smaller order.

\begin{theorem}[Iteration Complexity with Weighted Relative Error]\label{thm:itercomp} 
Define, for $\alpha \in (0, 1)$ and $\delta_0 > 0$, the (weighted) relative error class \begin{equation}\label{relerr-class} w\mathcal{R} := \left\{R(x,Y) \left(\|\nabla f(x)\|^{\alpha} + \delta_0\right), x \in \mathbb{R}^d\right\}; \quad R(x,Y): = \frac{\nabla F(x,Y) - \nabla f(x)}{\| \nabla f(x)\| + c_0}\end{equation} and its empirical-mean envelope $$H_m : = \sup_{x \in \mathbb{R}^d} \|R_{m}(x)\|\left( \|\nabla f(x) \|^{\alpha} + \delta_0\right).$$ Suppose, for some $\alpha >0$ and $\delta_0 >0$, the empirical-mean envelope of the weighted relative error class exhibits ``CLT-scaling,'' that is there exists $\sigma_0 < \infty$ such that,  \begin{equation}\label{emperrenv-clt} \mathbb{E}\left[\left(\sqrt{M_k}\, H_{M_k}\right)^{1/\alpha} \, \vert \, \mathcal{F}_{k-1}\right] \leq \sigma_0^{1/\alpha}; \quad M_k \in \mathcal{F}_{k-1}.\end{equation}  Then, for $k \geq 2$ and $\kappa_2:= 2(1-\alpha)\left(\frac{2\alpha}{1-\alpha}\right)^{\alpha} \sigma_0,$ we have that \begin{equation}\label{geniter} \mathbb{E}\left[ \| \nabla f (X_k) \| \, \vert \, \mathcal{F}_{k-1} \right] \leq  \frac{2c_2}{M_k^p} + \frac{2c_0\sigma_0/\delta_0}{\sqrt{M_k}} +  \left(\frac{\kappa_2}{\sqrt{M_k}}\right)^{1/\alpha}. \end{equation} Choose the sample-size sequence $\{M_k, k \geq 1\}$ according to~\eqref{sampsizeseqchoice}, and error tolerance according to~\eqref{errtolchoice} with $p=1/2$. Then for $k \geq 2,$ 
\begin{align} \label{itercompalt}
    \mathbb{E} \left[\|\nabla f(X_k)\|\right] 
    \leq \bigg(\frac{1}{\sqrt{c_{1}}}\bigg)^{k-1} \bigg( \frac{2(c_2 + c_0\sigma_0/\delta_0)}{\sqrt{m_1}} + \left(\frac{\kappa_2}{\sqrt{m_1}}\right)^{1/\alpha} \bigg).
\end{align}
\end{theorem}

Its proof will be given in Section \ref{ProofAltThm} of the Appendix. The power error class in Section \ref{PECEx} provides an example illustrating its applicability.

Analogous to Theorem~\ref{thm:itercomp0}\eqref{epsiloncomp0}, Theorem~\ref{thm:itercomp} leads to the following simple corollary that provides an $\epsilon$-complexity result.

\begin{corollary}[$\epsilon$-complexity with Weighted Relative Error]\label{epsiloncomp} Let the postulates  Theorem~\ref{thm:itercomp} along with the sample size choice in~\eqref{sampsizeseqchoice} and the error tolerance choice in~\eqref{errtolchoice} with $p=1/2$. Then for all $k \geq K(\epsilon)$ where $$K(\epsilon):= 1 + \frac{1}{\log \sqrt{c_1}}\left(\log \frac{1}{\epsilon} + \log \left(\frac{2(c_2 + c_0\sigma_0/\delta_0)}{\sqrt{m_1}} + \left(\frac{\kappa_2}{\sqrt{m_1}}\right)^{1/\alpha}\right)\right),$$ we have that $$\mathbb{E}\left[ \| \nabla f (X_k) \| \, \vert \, \mathcal{F}_{k-1} \right] \leq \epsilon, \quad \epsilon > 0.$$
\end{corollary}

\subsection{Oracle Complexity}\label{sec:oraccomp} 

To understand oracle complexity, recall that unlike SGD, RA algorithms may expend a varying number of oracle calls during each outer iteration, that is, between subsample updates. Specifically, the total number of oracle calls expended by RA after $k$ iterations is: \begin{equation} \label{totalwork} W_k : = \sum_{j=1}^k N_j \times M_j,\end{equation} where $N_j$ is the number of solver steps (or inner iterations), and $M_j$ is the adapted sample-size during the $j$-th iteration. Notice that the expression for the total work appearing in~\eqref{totalwork} assumes a ``linear sum'' structure; such structure is not always valid but holds for many common optimization solvers such as gradient descent~\cite[Chapter 2]{2004nes}, accelerated gradient descent~\cite[Chapter 2]{2004nes}, momentum method~\cite[Section 7]{bottou2018optimization}, and line search with L-BFGS~\cite[Algorithm 7.5]{nocedalbook}, since these methods require only gradient evaluations at each inner iteration.

The oracle complexity attained by RA depends on the speed of the solver employed to execute the inner iterations. Accordingly, we present our oracle complexity results as a function of whether the solver used in the inner iterations is \emph{linearly convergent} or \emph{sub-linearly convergent.} --- see Definition~\ref{defn:solverspeed}. Theorem~\ref{thm:expworkcomlin} that follows is intended to address problems having strongly convex sample-path problems while Theorem~\ref{thm:expworkcomsublin} addresses problems having (merely) convex or nonconvex sample-path problems.

\begin{theorem}[Expected Oracle Complexity for Linear Solvers]\label{thm:expworkcomlin} Suppose the Method $\mathfrak{M}$ within RA exhibits $(\rho^{\tiny{\mathfrak{M}}}_{M_k},C^{\tiny{\mathfrak{M}}}_{M_k})$-linear convergence (see Definition~\ref{defn:solverspeed}) when executed on  the sample path problem $(S_{M_k})$ starting at initial point $X_{k-1}$, where $0 < \rho_{M_k}^{\scaleto{\mathfrak{M}}{4pt}} \leq \rho_{\scaleto{\mathfrak{M}}{4pt}} < 1$ and $C_{M_k}^{\scaleto{\mathfrak{M}}{4pt}} \leq C_{\scaleto{\mathfrak{M}}{4pt}}$, and $\rho_{\scaleto{\mathfrak{M}}{4pt}}< 1, C_{\scaleto{\mathfrak{M}}{4pt}} < \infty$ are some constants. Let the sample size sequence $\{M_k, k \geq 1\}$ and the error tolerance sequence $\{\epsilon_k, k \geq 1\}$ be chosen so that \begin{align}\label{sampsizeerrtolchoice} M_k &:= C_{1,k}M_{k-1}, & C_{1,k} = c_1, 1 < c_1 < \infty; \\ \epsilon_k &:= \frac{C_{2,k}}{\sqrt{M_k}}, &  C_{2,k} \in [\underline{c}_2,c_2], 0 < \underline{c}_2 \leq c_2 < \infty, \label{eq:errortoler}\\ \tag{C.3} \end{align} where $C_{1,k}, C_{2,k} \in \mathcal{F}_{k-1}$, and $c_1, \bar{c}_1, c_2, \underline{c}_2$ are (non-random) constants.

If the postulates of Theorem~\ref{thm:itercomp0} hold, then for each $k \geq 1$,
\begin{align}\label{implication:oraccomp}
    \E[W_k] = O\left(\frac{A^2 m_1}{(c_1 - 1)\epsilon^2\log(1/\rho_{\scaleto{\mathfrak{M}}{4pt}})}\right) \,\,\implies\,\, \mathbb{E}\left[\nabla f(X_k)\right] \leq \epsilon,
\end{align}
with \begin{align}\label{eq:Athm3}
        A = \frac{\sqrt{c_1}(c_2 + c_0\sigma + \sigma\ell_{2,\infty} )}{\sqrt{m_1}}.
    \end{align} 
    
    If the postulates of Theorem~\ref{thm:itercomp} hold and ~\eqref{emperrenv-clt0} is satisfied, then~\eqref{implication:oraccomp} holds with  
    \begin{align}\label{eq:Athm4} A =  \sqrt{c_1}\left(\frac{2(c_2 + c_0\sigma_0/\delta_0)}{\sqrt{m_1}}\right)  + \left(\frac{\kappa_2}{\sqrt{m_1}}\right)^{1/\alpha}.
    \end{align} 
\end{theorem}

\begin{proof}{Proof.}
    From the definition of linear convergence (Definition~\ref{defn:solverspeed}), we see that for any
    \begin{align*}
        N_k \geq \frac{\log{\left(C_{\scaleto{\mathfrak{M}}{4pt}}\|\nabla f_{M_k}(X_{k-1})\|/\epsilon_k\right)}}{\log(1/\rho_{\scaleto{\mathfrak{M}}{4pt}})}
    \end{align*}
    the termination condition $\|\nabla f_{M_k} (X_k)\| \leq \epsilon_k$ is satisfied. 
    We can choose $N_k$ to be the smallest number such that the condition $\|\nabla f_{M_k} (X_k)\| \leq \epsilon_k$ is satisfied. Therefore, 
    \begin{align*}
        N_k = \ceil*{\frac{\log{\left(C_{\scaleto{\mathfrak{M}}{4pt}}\|\nabla f_{M_k}(X_{k-1})\|/\epsilon_k\right)}}{\log(1/\rho_{\scaleto{\mathfrak{M}}{4pt}})}} \leq 1 +  \frac{\log{\left(C_{\scaleto{\mathfrak{M}}{4pt}}\|\nabla f_{M_k}(X_{k-1})\|/\epsilon_k\right)}}{\log(1/\rho_{\scaleto{\mathfrak{M}}{4pt}})}
    \end{align*}
    Taking expectation, using~\eqref{eq:Gm}, and applying Jensen's inequality on a concave function, we have,
    \begin{align*}
        \E[N_k]&\leq 1 + \frac{\log(C_{\scaleto{\mathfrak{M}}{4pt}}) + \log(\E[\|\nabla f_{M_k}(X_{k-1})\|/\epsilon_k])}{\log(1/\rho_{\scaleto{\mathfrak{M}}{4pt}})} \\
        &\leq 1 + \frac{\log(C_{\scaleto{\mathfrak{M}}{4pt}}) + \log(\E[\|\nabla f_{M_k}(X_{k-1}) - \nabla f(X_{k-1})\| + \|\nabla f(X_{k-1})\|/\epsilon_k])}{\log(1/\rho_{\scaleto{\mathfrak{M}}{4pt}})} \\
        &\leq 1 + \frac{\log(C_{\scaleto{\mathfrak{M}}{4pt}}) + \log(\E[\|G_{M_k}\|(\|\nabla f(X_{k-1})\| + c_0) + \|\nabla f(X_{k-1})\|/\epsilon_k])}{\log(1/\rho_{\scaleto{\mathfrak{M}}{4pt}})} \\
        &\leq 1 + \frac{\log(C_{\scaleto{\mathfrak{M}}{4pt}})}{\log(1/\rho_{\scaleto{\mathfrak{M}}{4pt}})} + \frac{\log \left(\sigma \E[(\|\nabla f(X_{k-1})\| + c_0)] + \sqrt{M_k}\E[\|\nabla f(X_{k-1})\|]\right) - \log{\underline{c_2}}}{\log(1/\rho_{\scaleto{\mathfrak{M}}{4pt}})} 
    \end{align*} 
    where the last inequality is due to the assumption~\eqref{emperrenv-clt0} and \eqref{eq:errortoler}.
    Now from Theorem~\ref{thm:itercomp0} and~\ref{thm:itercomp}, we have that 
    \begin{align*}
        \E[\|\nabla f(X_{k-1})\|] \leq C\left(\frac{1}{\sqrt{c_1}}\right)^{k-1}. 
    \end{align*}
    Substituting this in the above inequality, we get,
    \begin{equation*}
        \E[N_k]\leq 1 + \frac{\log(C_{\scaleto{\mathfrak{M}}{4pt}})}{\log(1/\rho_{\scaleto{\mathfrak{M}}{4pt}})} + \frac{\log\left(\sigma(Cc_1^{-(k-1/2)} + c_0) + \sqrt{m_1}C\right) - \log{\underline{c_2}}}{\log(1/\rho_{\scaleto{\mathfrak{M}}{4pt}})} =:
 \frac{C_u}{\log(1/\rho_{\scaleto{\mathfrak{M}}{4pt}})}
    \end{equation*}
    where 
    \begin{align*}
        C_u = \log\left(\frac{C_{\scaleto{\mathfrak{M}}{4pt}}\left(\sigma(Cc_1^{-(k-1/2)} + c_0) + \sqrt{m_1}C\right)}{\rho_{\scaleto{\mathfrak{M}}{4pt}}\underline{c_2}}\right).
    \end{align*}
    The expected work complexity to achieve $\E[\|\nabla f(X_k)\|] \leq \epsilon$ is then
    \begin{align}\label{exporaccompexp}
        \E[W_{K(\epsilon)}] &=\sum_{j=1}^{{K(\epsilon)}} M_j \E[N_j] \leq \sum_{j=1}^{{K(\epsilon)}} m_1\bar{c}_1^{j-1}\frac{C_u}{\log(1/\rho_{\scaleto{\mathfrak{M}}{4pt}})} \leq m_1\frac{C_u}{\log(1/\rho_{\scaleto{\mathfrak{M}}{4pt}})} \frac{c_1^{K(\epsilon)} - 1}{c_1 - 1}.
    \end{align}
    
    Definition~\eqref{epsiloncomp0} and Corollary~\ref{epsiloncomp} imply 
    $K(\epsilon) = \frac{\log\left(\frac{A}{\epsilon}\right)}{\log \sqrt{c_1}}$ and
$c_1^{K(\epsilon)} = (\frac{A}{\epsilon})^2.$
Plugging in~\eqref{exporaccompexp}, we get 
\begin{align*}
     \E[W_{K(\epsilon)}] &\leq  \frac{m_1}{c_1 - 1} \frac{C_u}{\log(1/\rho_{\scaleto{\mathfrak{M}}{4pt}})}c_1^{K(\epsilon)}
     \leq  \frac{m_1}{c_1 - 1} \frac{C_u}{\log(1/\rho_{\scaleto{\mathfrak{M}}{4pt}})} \left(\frac{A}{\epsilon}\right)^2.  \blacksquare 
\end{align*}

\end{proof}

\begin{remark}
It is important that the result obtained through Theorem~\ref{thm:expworkcomlin} is predicated on the existence of a linearly convergent solver. Such a convergence rate is usually exhibited by solvers when executed on strongly convex functions, or on nonconvex functions when a PL-type condition~\cite{1963pol} holds. As such, we are unaware of schemes which provide such non-asymptotic guarantees otherwise.
\end{remark}

Observe that Theorem~\ref{thm:expworkcomlin} signals the optimal $O(\epsilon^{-2})$ oracle complexity, with the constant $A$ involving the rate of increase of sample size $c_1$,  the inherent ``variability'' $\sigma$ of the stochastic oracle, and the extent of the excursions $\ell_{2,\infty}$ of RA's iterates. Theorem~\ref{thm:expworkcomlin} assumes that $C_{1,k}=c_1$, that is, the sample sizes are increased at a fixed rate $c_1;$ it can be shown that allowing the rate of increase of the sample sizes to be in an interval introduces a slight inflation in the expected oracle complexity rate.    

Next, we prove a result analogous to Theorem~\ref{thm:expworkcomlin} but when a sub-linear solver is used when solving the sample-path problems.

\begin{theorem}[Expected Oracle Complexity for Sublinear Solvers]\label{thm:expworkcomsublin}
Suppose the Method $\mathfrak{M}$ within RA exhibits $(1/2,C^{\tiny{\mathfrak{M}}}_{M_k})$-sublinear convergence (see Definition~\ref{defn:solverspeed}) when executed on  the sample path problem $(S_{M_k})$ starting at initial point $X_{k-1}$, where $0 < C_{M_k}^{\scaleto{\mathfrak{M}}{4pt}} \leq C_{\scaleto{\mathfrak{M}}{4pt}}$, and $ C_{\scaleto{\mathfrak{M}}{4pt}} < \infty$ is some constant. Let the sample size sequence$\{M_k, k \geq 1\}$ and the error tolerance sequence $\epsilon_k, k \geq 1\}$ be chosen according to~\eqref{sampsizeerrtolchoice}, and suppose further that there exists $C_A< \infty$ such that \begin{equation}\label{fniterbd} \sup_k \, \mathbb{E}[f(X_k)] < \frac{C_A}{2}; \quad \, \mathbb{E}[\inf_{x \in \mathbb{R}^d} F(x,Y)] > -\frac{C_A}{2}.\tag{\mbox{Bd-Fn}} \end{equation} If the postulates of Theorem~\ref{thm:itercomp0} or Theorem~\ref{thm:itercomp} hold and~\eqref{emperrenv-clt0} is satisfied, then, for each $k \geq 1$,
\begin{align}\label{implication:oraccompsublin}
    \E[W_k] = O\left(\frac{A^4 m_1 C_{\scaleto{\mathfrak{M}}{4pt}}C_{A}}{(c_1^2 - 1)\underline{c_2}^{2}\epsilon^4}\right) \,\,\implies\,\, \mathbb{E}\left[\nabla f(X_k)\right] \leq \epsilon,
\end{align} 
where $A$ is as in~\eqref{eq:Athm3} or~\eqref{eq:Athm4} (according to whether Theorem~\ref{thm:itercomp0} or Theorem~\ref{thm:itercomp} hold),  and 
\begin{equation}\label{eq:cu}
    C_u = \frac{C_{\scaleto{\mathfrak{M}}{4pt}}\left(\sigma(Cc_1^{-(k-1/2)} + c_0) + \sqrt{m_1}C\right)}{\underline{c_2}}. 
\end{equation}

\end{theorem}

\begin{proof}{Proof.}
From the definition of sublinear convergence (Definition~\ref{defn:solverspeed}), we see that for any 
\begin{equation*}
    N_k \geq \frac{C_{\scaleto{\mathfrak{M}}{4pt}} (f_{M_k}(X_{k-1}) - f_{M_k}(X_k))}{\epsilon_k^2}
\end{equation*}
the termination condition $\|\nabla f_{M_k} (X_k)\| \leq \epsilon_k$ is satisfied. We can choose $N_k$ to be the smallest number such that $\|\nabla f_{M_k} (X_k)\| \leq \epsilon_k$ condition is satisfied. Therefore, 
\begin{align*}
    N_k = \ceil*{\frac{C_{\scaleto{\mathfrak{M}}{4pt}} (f_{M_k}(X_{k-1}) - f_{M_k}(X_k))}{\epsilon_k^2}} \leq  \frac{2C_{\scaleto{\mathfrak{M}}{4pt}} (f_{M_k}(X_{k-1}) - f_{M_k}(X_k))}{\epsilon_k^2}
\end{align*}
Taking expectation, we have,
\begin{align*}
    \E[N_k]&\leq \E\left[\frac{2C_{\scaleto{\mathfrak{M}}{4pt}} (f_{M_k}(X_{k-1}) - f_{M_k}(X_k))}{\epsilon_k}\right] \\
    &\leq \frac{2M_k C_{\scaleto{\mathfrak{M}}{4pt}} (\E[f_{M_k}(X_{k-1}) - f_{M_k}(X_{k})])}{\underline{c_2}^{2}} 
\end{align*}
Now using \eqref{fniterbd}, we can show 
\begin{align*}
    \E[f_{M_k}(X_{k-1}) - f_{M_k}(X_{k})] &\leq \E[f(X_{k-1})] - \E[\inf_{x \in \mathbb{R}^d}f_{M_k}(x)] \leq C_{A}. 
\end{align*}
Therefore, we have, 
\begin{align*}
     \E[N_k]&\leq  \frac{2M_k C_{\scaleto{\mathfrak{M}}{4pt}}C_{A}}{\underline{c_2}^{2}} 
\end{align*}

  Now from Theorem~\ref{thm:itercomp0} and~\ref{thm:itercomp}, we have that 
\begin{align*}
    \E[\|\nabla f(X_{k-1})\|] \leq C\left(\frac{1}{\sqrt{c_1}}\right)^{k-1}. 
\end{align*}
Following the same steps employed in Theorem~\ref{thm:expworkcomlin}, we obtain
\begin{align*}
    \E[W_{K(\epsilon)}] &=\sum_{j=1}^{{K(\epsilon)}} M_j \E[N_j] \leq  \sum_{j=1}^{{K(\epsilon)}} m_1^2c_1^{2(j-1)}\frac{2 C_{\scaleto{\mathfrak{M}}{4pt}}C_{A}}{\underline{c_2}^{2}}
    \leq \frac{2 m_1^2C_{\scaleto{\mathfrak{M}}{4pt}}C_{A}}{\underline{c_2}^{2}}\frac{c_1^{2K(\epsilon)} - 1}{c_1^2 - 1}.
\end{align*}
From~\eqref{epsiloncomp0} and Corollary~\ref{epsiloncomp}, we have 
$K(\epsilon) = \frac{\log\left(\frac{A}{\epsilon}\right)}{\log \sqrt{c_1}}$ and 
substituting this into the above inequality gives the desired result. $\blacksquare$
\end{proof}

The derived complexities in Theorem~\ref{thm:expworkcomlin} for linear solvers and Theorem~\ref{thm:expworkcomsublin} for sublinear solvers should be loosely interpreted as $O(\epsilon^{-2})$ and $O(\epsilon^{-4}),$ respectively. The former complexity matches the best achievable rate for stochastic optimization on smooth strongly convex functions, while the latter rate matches the  $O(\epsilon^{-4})$ lower bound for smooth nonconvex stochastic optimization~\cite{arjevani2019lower}. 

Also, Theorem~\ref{thm:expworkcomsublin} relies on the technical condition in~\eqref{fniterbd}. To understand this condition, notice that we have assumed little about the nature of the solutions obtained by the solvers used within RA on successive iterations. So, presumably, since $f$ can be nonconvex, the solver in use and the problem landscape might collude in a way that RA's trajectory oscillates between first-order critical point estimators in regions with a high objective and those with a low objective. The condition in~\eqref{fniterbd} prevents such oscillation from being extreme.  

Recall that the $k$-th (outer) iteration of RA is terminated when the norm $\| \nabla f_{M_k}(X_k)\|$ of the subsampled gradient at the incumbent iterate falls below the threshold $\epsilon_k$, that is, if $X_{k-1,t}$ is the $t$-th iterate obtained from the employed solver during the $k$-th iteration, then the number of inner steps $$N_k := \min\left\{t: \|\nabla f_{M_k}(X_{k,t})\| \leq \epsilon_k\right\}.$$  Theorem~\ref{thm:expworkcomlin} and Theorem~\ref{thm:expworkcomsublin} consider the context where the tolerance $\epsilon_k$ is chosen as $\epsilon_k = C_{2,k}/\sqrt{M}_k,$ implying that \begin{equation}\label{stopgeneric}N_k:= \min\{t: \|\nabla f_{M_K}(X_{k,t})\| \leq \frac{C_{2,k}}{\sqrt{M_k}}\}.\end{equation}

An alternative to the stopping criterion in~\eqref{stopgeneric}, and one that is consistent with practice within deterministic solvers, is to terminate the inner iterations of the $k$-th outer iteration when the relative error drops below a threshold, that is, \begin{equation}\label{errtolrelerr}N_k := \min\left\{t: \|\nabla f_{M_k}(X_{k,t})\| \leq \Gamma_k\,\|\nabla f_{M_k}(X_{k-1})\|\right\},\end{equation} where $\Gamma_k \in \mathcal{F}_{k-1}.$ 

To characterize the expected oracle complexity rate when the stopping criterion is chosen according to the ``practical condition'' in~\eqref{errtolrelerr}, we need the following alternate model of progress for sublinear solvers. A solver is said to exhibit \emph{$C$-sublinear convergence} ($C < \infty$) if the gradient norm at iterate $x_n$ obtained after executing $n$ steps of the algorithm satisfies \begin{equation}\label{sublinconv-alt}\| \nabla f(x_n) \|  \leq C \, \| \nabla f(x_0)\|\, n^{-1/2}.\end{equation} The model of progress in~\eqref{sublinconv-alt} is possibly more typical of behavior in practice than that in~\eqref{sublinconv} but there appears to be no straightforward justification such as in~\eqref{graddescentconvrate}. Using the alternate notion of sublinear convergence in~\eqref{sublinconv-alt}, we can make the following remark about the oracle complexity when using the practical termination condition in~\eqref{errtolrelerr}. 
\begin{remark}\label{rm:workcomp_alt}
Suppose we set \begin{equation*}\label{eq:relterm}
    \epsilon_k := \Gamma_k \|\nabla f_{M_k}(X_{k-1})\|,
\end{equation*}
where $\Gamma_k \geq \gamma > 0, \Gamma_k \in \mathcal{F}_{k-1}$. Let the postulates of Theorem~\ref{thm:itercomp0} or Theorem~\ref{thm:itercomp} hold and let \eqref{emperrenv-clt0} be satisfied. If a linear solver or a sublinear solver satisfying \eqref{sublinconv-alt} is used within RA, then, for each $k\geq1$, 
\begin{align*}
    \E[W_k] \geq O(\epsilon^{-2}) \implies \E[\nabla f(X_k)] \leq \epsilon. 
\end{align*}
\end{remark}
Since the above statement can be proved using the same steps as in Theorems~\ref{thm:itercomp0},~\ref{thm:itercomp} and \ref{thm:expworkcomlin} we do not include a proof.

\section{EXAMPLE PROBLEM CLASSES}\label{sec:examples} It should be clear that the uniform law~\eqref{ass:RULLN} on the relative error class (or the sub-linear entropy growth condition in~\eqref{ass:sublinentgrowth}, as detailed in the Appendix) and the bounded gradients assumption appearing in~\eqref{graditerbd} essentially ensure in-probability, almost-sure, and $L_1$ convergence of RA's iterates. How stringent are~\eqref{ass:RULLN} and~\eqref{graditerbd}? To further the reader's intuition in support of our claim that~\eqref{ass:RULLN} and~\eqref{graditerbd}  are widely satisfied, we now describe a number of example problem classes, in each case discussing whether~\eqref{ass:RULLN} and~\eqref{graditerbd} are satisfied. While we do not prove, the ``CLT-scaling'' assumption employed in Theorems~\ref{thm:itercomp0}--\ref{thm:expworkcomlin} also hold for each of these classes under further mild assumptions.

\subsection{Example (Uniform Error Class).}\label{ex:uniferr} Suppose that the error in the gradient norm is uniformly convergent, that is, \begin{equation}\label{ulln} \sup_{x \in \mathcal{X} \subseteq \mathbb{R}^d} \| \nabla f_m(x) - \nabla f(x) \| \as 0 \mbox{ as } m \to \infty, \tag{ULLN}\end{equation} and that the envelope $$H = H(Y) := \sup\{\|\nabla F(x,Y) - \nabla f(x) \|: x \in \mathcal{X}\}$$ has finite second moment, that is, \begin{equation} \label{absfinsecmom} \mathbb{E}[H^2] < \infty.\end{equation} Variations of such a uniform error assumption are standard in SAA and empirical risk minimization literature, e.g., see~\cite[pp. 157,174]{2009shadenrus} and~\cite{2019mokozdjad,2006barjormca,2011botbos,2018meibaimon}.) Under~\eqref{ulln} and~\eqref{absfinsecmom}, notice that \begin{equation}\label{uniferrineq}\sup_{x \in \mathcal{X}}\|R_m(x)\| = \sup_{x \in \mathcal{X}}\frac{\|\nabla f_m(x) - \nabla f(x)\|}{\|\nabla f(x)\| + c_0} \leq \sup_{x \in \mathcal{X}}\frac{\|\nabla f_m(x) - \nabla f(x)\|}{c_0} \as 0,\end{equation} implying that~\eqref{ass:RULLN} is satisfied, and that strong consistency is guaranteed if RA is implemented on this class with $M_k \as \infty$ and $\epsilon_k \as 0$. 

Let's next demonstrate that the condition in~\eqref{graditerbd} needed for $L_1$ consistency of RA's iterates is satisfied as well. First, note that due to the assumption in~\eqref{ulln}, the $\delta$-covering number of the absolute error class $$\mathcal{A} := \{\nabla F(x,Y) - \nabla f(x), x \in \mathcal{X}\}$$ exhibits sub-linear empirical entropy growth, that is, \begin{equation}\label{abserrgrowth} \frac{1}{M_k} H_1(\delta,\mathcal{A},P_{M_k}) \inP 0 \mbox{ as } M_k \as \infty.\end{equation} Due to~\eqref{abserrgrowth}, we see that the postulates of Theorem 2.4.3~\cite[pp. 123]{1996vanwel} hold, and we have \begin{equation}\label{uniferrol1bd} \mathbb{E}\left[ \sup_{x \in \mathcal{X}} \, \| \nabla f_{M_k}(x) - \nabla f(x) \| \, \vert \, \mathcal{F}_{k-1} \right] \as 0 \mbox{ as } k \to \infty. \end{equation}  Recall that the iterates $\{X_k, k \geq 1\}$ generated by RA are such that \begin{equation}\label{termcritagain} \| \nabla f_{M_k}(X_k) \| \leq \epsilon_k.\end{equation} Using~\eqref{uniferrol1bd} and~\eqref{termcritagain}, and since $\epsilon_k \to 0$, we conclude that the key condition in~\eqref{graditerbd} holds.    
$\blacksquare$

Indeed, while~\eqref{ulln} is a common assumption, it is too stringent and is frequently violated by even standard stochastic optimization settings. Let's illustrate this with a simple example. Suppose $\mathcal{X} = \mathbb{R}^d$, $b \in \mathbb{R}$, $a= \mathbb{E}[Z(Y)]$, and $Z(Y) > 0$ is some positive-valued random variable with finite expectation. Set \begin{equation}\label{counterex}F(x,Y) = \frac{1}{2}Z(Y)x^Tx + bx; \quad  f(x) = \frac{1}{2}ax^Tx + bx.\end{equation}
We then see that 
$$\nabla F(x,Y) = Z(Y)x + b; \quad \nabla f(x) = ax + b$$ and hence, if $\mbox{Var}(Z(Y)) \neq 0,$ then
$$\sup_{x \in \mathcal{R}^d}\| \nabla f_m(x) - \nabla f(x) \| = \left|\frac{1}{n} \sum_{j=1}^n Z(Y_j) - a\right |\sup_{x \in \mathcal{R}^d}\|x\| \not\to 0,$$ thus violating~\eqref{ulln}.

The main feature of the counter-example in~\eqref{counterex} responsible for~\eqref{ulln} violation is that the error in the sampled gradient is proportional to the true gradient, that is, \begin{equation}\label{counterexineq}\left \|\nabla F(x,Y) - \nabla f(x)\right\| = \left \| (Z(Y)-a) x \right \| \geq \left |\frac{Z(Y)-a}{a}\right |\left(\|\nabla f(x)\| - b\right).\end{equation} The inequality in~\eqref{counterexineq} implies that the absolute error in the gradient estimate for the counter-example in~\eqref{counterex} will increase proportionally with the norm of the true gradient, thereby violating~\eqref{ulln} if the true gradient is unbounded over $\mathcal{X}$. 

This issue of absolute error in the estimated gradient norm being proportional to the true gradient seems to be common in practice~\cite{bottou2018optimization}. Through the next example class, we demonstrate that in many simple settings, while a uniform law~\eqref{ulln} on the absolute error is violated, a uniform law~\eqref{ass:RULLN} on the relative error is satisfied. 

\subsection{Example (Multiplicative Error  Class).}\label{ex:multerrclass} Suppose the random function $\nabla F(\cdot,Y)$ is such that \begin{equation}\label{mulerrclass}\nabla F(x,Y) = A(Y)\nabla f(x) + B(Y), \quad x \in \mathbb{R}^d\end{equation} where $A(Y)$ is a $d \times d$ matrix of real-valued random variables and $B(Y)$ is an $\mathbb{R}^d$ valued random variable satisfying \begin{equation} \mathbb{E}[A(Y)] =J_d; \quad \mathbb{E}[B(Y)] =0, \end{equation} and $J_d$ is the $d \times d$ matrix of ones. Then, we see that \begin{equation} \label{multerrclass} R_m(x) = \frac{(\bar{A}_m - J_d)\nabla f(x) + \bar{B}_m}{\| \nabla f(x) \| + c_0}; \quad \bar{A}_m := \frac{1}{m} \sum_{j=1}^m A(Y_j); \quad \bar{B}_m := \frac{1}{m} \sum_{j=1}^m B(Y_j),\end{equation} implying that the empirical-mean envelope $G_m:= \sup_{x \in \mathcal{X}} \|R_m(x)\|$ of the relative error class satisfies \begin{equation}\label{mainmulineq} G_m \leq \left \|\bar{A}_m - J_d\right \| + \left \|\bar{B}_m \right \| \frac{1}{c_0}. \end{equation} Conclude from~\eqref{mainmulineq} that the uniform law in~\eqref{ass:RULLN} is satisfied. 

Let's next demonstrate that under further moment conditions, the condition in~\eqref{graditerbd} is satisfied as well. Suppose there exist deterministic constants $\ell_1,\ell_2< \infty$ such that \begin{equation}\label{l1} \mathbb{E}\left[ \left(\frac{1}{\sigma^*_k}\right)^{4} \, \vert \, \mathcal{F}_{k-1}\right] \leq \ell_1 < \infty; \quad \sigma^*_k := \min_{x \neq 0} \frac{\| \bar{A}_{M_k}x\|}{\|x\|} \end{equation} and \begin{equation}\label{l2}\mathbb{E}\left[ \|B(Y)\|^{4} \right] \leq \ell_2 < \infty,\end{equation} (Notice that $\sigma^*_k$ is the minimum singular value of the matrix $\bar{A}_{M_k}$, that is, the square root of the smallest eigenvalue of the symmetric matrix $\bar{A}_{M_k}\bar{A}_{M_k}^T$.) Suppose also that the error tolerance $\epsilon_k$ in RA is chosen so that \begin{equation}\label{l3} \mathbb{E}\left[\epsilon_k^4 \, \vert \, \mathcal{F}_{k-1}\right] \leq \ell_3. \end{equation} From~\eqref{l2}, we see that \begin{equation}\label{l1l2sampmean} \mathbb{E}\left[\left\|\bar{B}_{M_k}\right\|^4 \, \vert \, \mathcal{F}_{k-1}\right] \leq \ell_2.\end{equation} Next, we see from~\eqref{mulerrclass} that $X_k$ satisfies \begin{equation}\label{stoppingcond} \left\|\bar{A}_{M_k} \nabla f(X_k) + \bar{B}_{M_k}\right\| \leq \epsilon_k, \end{equation} implying that  \begin{equation}\label{gradineqfirst}\left\|\nabla f(X_k)\right\|  \leq \| \bar{A}_{M_k}^{-1}\|\left(\epsilon_k + \| \bar{B}_{M_k}\|\right).\end{equation} Also, \begin{equation}\label{norminv} \|\bar{A}_{M_k}^{-1}\| := \max_{x \neq 0} \frac{\|\bar{A}_{M_k}^{-1}x\|}{\|x\|} = \left( \min_{x \neq 0} \frac{\|\bar{A}_{M_k}x\|}{\|x\|}\right)^{-1} = \frac{1}{\sigma^*_k}.\end{equation} Hence \begin{align} \label{supl2bdcalc} \MoveEqLeft \sup_{k} \left\{\mathbb{E}\left[\left\|\nabla f(X_k)\right\|^2 \, \vert \, \mathcal{F}_{k-1} \right]\right\} \nonumber  & \nonumber \\ & \leq \sup_{k} \left\{2\,\mathbb{E}\left[ \left(\frac{1}{ \sigma^*_{k}}\right)^2\epsilon_k^2 \, \vert \, \mathcal{F}_{k-1}\right] + 2\,\mathbb{E}\left[\left(\frac{1}{\sigma^*_{k}}\right)^2 \|\bar{B}_{M_k} \|^2\, \vert \, \mathcal{F}_{k-1}\right] \right\}\nonumber \\ & \leq 2\sup_{k} \left\{ \left(\mathbb{E}\left[ \left(\frac{1}{ \sigma^*_{k}}\right)^4 \, \vert \, \mathcal{F}_{k-1}\right] \right)^{1/2}\left(\left(\mathbb{E}\left[\epsilon_k^4 \, \vert \, \mathcal{F}_{k-1}\right]\right)^{1/2} + \left(\mathbb{E}\left[\left\|\bar{B}_{M_k}\right\|^4 \, \vert \, \mathcal{F}_{k-1}\right] \right)^{1/2}\right)\right\}.\end{align} Using~\eqref{l1},~\eqref{l2} and~\eqref{l3} in~\eqref{supl2bdcalc}, we see that the condition in~\eqref{graditerbd} is satisfied.$\blacksquare$

\begin{remark} Example~\ref{ex:multerrclass} is a special case of a larger and more useful class of random functions having random coefficients that are dependent on $x$, that is, $$\nabla F(x,Y) = A(x,Y)\nabla f(x) + B(x,Y).$$ An analysis similar to what we have outlined in Example~\ref{ex:multerrclass} is possible for this class but with stipulations on $\{A(x,Y), x \in \mathbb{R}^d\}$ and $\{B(x,Y), x \in \mathbb{R}^d\}$ replacing our current stipulations on $A(Y)$ and $B(Y)$.
\end{remark}



\subsection{Example (Linear and Polynomial Regression with Least Squares).}\label{ex:regression} Let's consider the classical regression problem of statistics where we wish to estimate the parameters $x^* \in \mathbb{R}^d$ of an assumed linear regression model \begin{equation}\label{regressionmodel}T(Y) = Qx^* + \epsilon(Y),\end{equation} where $T \in \mathbb{R}$ is the ``observed'' dependent variable, $Q = (Q^{(1)},Q^{(2)}, \ldots,Q^{(d)})$ is a $1 \times d$ random vector of covariates, $x^*$ is a $d \times 1$ vector of true parameters to be estimated using data. Assume:  \begin{equation}\label{ass:regression} \mathbb{E}[\epsilon(Y)]=0; \quad \mbox{Corr}(Q^{(j)},\epsilon(Y)) = 0 \,\, \forall \, j; \quad \mathbb{E}[Q^TQ] < \infty; \quad \mathbb{E}[Q^TT] < \infty; \mbox{ and } \left(\mathbb{E}[Q^TQ]\right)^{-1} \mbox{ exists. }  \end{equation} 

The question of estimating $x^*$ can be posed as the following unconstrained stochastic optimization problem: \begin{align}\label{regression} \mbox{minimize: } & f(x) := \frac{1}{2}\mathbb{E}\left[ (T - Qx)^2 \right] \nonumber \\ \mbox{ s.t. } & x \in \mathbb{R}^d.\end{align} Using~\eqref{ass:regression}, it can be shown that $x^*$ is the unique solution of the problem in~\eqref{regression}.

Suppose now that we observe data $(Q_j,T_j), i = 1,2,\ldots, n$, where $Q_j,T_j, j =1,2,\ldots,n$ are iid copies of $(Q,T)$ and $n > d$. Let $\tilde{Q}_n$ represent the $n \times d$ matrix obtained by stacking $Q_j, j=1,2,\ldots,n$ and let $\tilde{T}_n$ be the $n \times 1$ matrix obtained by stacking $T_j, j =1,2,\ldots,n$. (The ``tall'' $n \times d$ matrix $\tilde{Q}_n$ obtained by stacking $Q_i, i=1,2,\ldots,n$ is traditionally called the ``data matrix''~\cite{2009tro}.) The sample-path  problem corresponding to~\eqref{regression} is \begin{align}\label{sampregression} \mbox{minimize: } & f_n(x) := \frac{1}{2}\frac{1}{n} \sum_{j=1}^n (T_j - Q_jx)^2 \nonumber \\ \mbox{ s.t. } & x \in \mathbb{R}^d. \nonumber \\ \end{align} Now observe that \begin{align} \label{truegradls} \nabla f(x) &= \mathbb{E}\left[Q^TQ\right]x - \mathbb{E}\left[Q^TT\right];\nonumber \\ \nabla f_n(x) &= \frac{1}{n}\sum_{j=1}^n \left(Q_j^TQ_jx - Q_j^TT_j\right)= \frac{1}{n} \left(\tilde{Q}_n^T\tilde{Q}_nx - \tilde{Q}_n^T\tilde{T}_n\right). \nonumber \\ \end{align} Notice also that under the further assumption that the data matrix $\tilde{Q}_n$ has rank $d$, we get the classical least squares estimator $X_n^*$ for $x^*$:  \begin{equation}\label{samplepathsoln} X_n^* = \left(\tilde{Q}_n^T\tilde{Q}_n^{-1}\right)^{-1} \tilde{Q}_n^T\tilde{T}_n.\end{equation}

From~\eqref{truegradls}, we see that the relative error class $$\mathcal{R} := \{\left(Q^TQ - \mathbb{E}\left[Q^TQ\right]\right)x - \left(Q^TT - \mathbb{E}\left[Q^TT\right]\right), x \in \mathbb{R}^d\},$$ and the relative error $R_m(x)$ for any $x \in \mathbb{R}^d$ satisfies \begin{align}\label{relerrorregression} \left\| R_m(x) \right\| &= \frac{\left\| \bar{A}_mx -\bar{B}_m \right\|}{\left\|\mathbb{E}\left[Q^TQ\right]x - \mathbb{E}\left[Q^TT\right]\right\| + c_0} \nonumber \\ & \leq \|\bar{A}_m\|\frac{\|x\|}{\left\|\mathbb{E}\left[Q^TQ\right]x - \mathbb{E}\left[Q^TT\right]\right\| + c_0} + \frac{\|\bar{B}_m\|}{c_0},\end{align} where \begin{equation} \bar{A}_m := \mathbb{E}\left[Q^TQ\right] - \frac{1}{m}\sum_{j=1}^n Q_j^TQ_j; \quad \bar{B}_m := \mathbb{E}\left[Q^TT\right] - \frac{1}{m}\sum_{j=1}^n Q_j^TT_j.\end{equation} Due to the assumption in~\eqref{ass:regression}, \begin{equation}\label{regconsts} \bar{A}_m \as 0; \quad \bar{B}_m \as 0,\end{equation} and since $\mathbb{E}[Q^TQ]$ is nonsingular, we can show after some algebra that \begin{equation}\label{supfrac} \sup_{x \in \mathbb{R}^d}\, \frac{\|x\|}{\left\|\mathbb{E}\left[Q^TQ\right]x - \mathbb{E}\left[Q^TT\right]\right\| + c_0} < \infty.\end{equation} Using~\eqref{supfrac} and~\eqref{regconsts} in~\eqref{relerrorregression}, we see that the empirical-mean envelope $G_m := \sup_{x \in \mathbb{R}^d} \|R_m(x)\| \as 0$ and hence the uniform law in~\eqref{ass:RULLN} is satisfied assuring almost sure convergence of RA when implemented on this class.  

Let's now demonstrate that under further conditions, the bounded gradients assumption in~\eqref{graditerbd} is also satisfied, ensuring $L_1$ convergence of RA. Suppose that the following technical condition on the observed data holds for each $i = 1,2, \ldots, d$: \begin{equation}\label{techcondregression} \textcolor{blue}{\mathbb{E}\left[ \left( \left| M_k^{-1}\sum_{j=1}^{M_k} Q_j^{(i)}T_j  \right| - \epsilon_k \right)^{-4}\right] < \infty.}\end{equation} \textcolor{blue}{(To see that the technical condition in~\eqref{techcondregression} is not too stringent, consider the hypothetical context where the sample-path problem is solved to infinite precision, that is, $\varepsilon_k = 0$.)}  
After some algebra, we can write the sample gradient $\nabla f_n$ in terms of the true gradient $\nabla f$ as follows:

\begin{equation}\label{sampgradtograd} \nabla f_n(x) = \Lambda_n(x) \nabla f(x) + 2\left(\Lambda_n \mathbb{E}\left[Q^T T\right] - \tilde{Q}_n^T \tilde{T}_n \right),\end{equation} where $\Lambda_n$ is a $d \times d$ diagonal matrix given by \begin{equation}\label{Lambda} \Lambda_n(x): =\begin{bmatrix}
    \frac{n^{-1}\sum_{j=1}^n Q_j^{(1)}Q_jx}{\mathbb{E}\left[Q^{(1)}Q\right]x} & & \\
     & \ddots &  \\
    & & \frac{n^{-1}\sum_{j=1}^n Q_j^{(d)}Q_jx}{\mathbb{E}\left[Q^{(d)}Q\right]x}
  \end{bmatrix}. \end{equation} From~\eqref{truegradls} and since the solution $X_k$ obtained after $k$ iterations of RA satisfies \begin{equation}\label{samplepathsoln} \| \nabla f_{M_k}(X_k) \| \leq \epsilon_k, \end{equation} we see that \begin{equation}\label{truegradatXk} \left\|\nabla f(X_k)\right\| \leq \left\| \Lambda^{-1}_{M_k}(X_k) \right\| \left( \epsilon_k + 2\left\| \tilde{Q}_{M_k}^T\tilde{T}_{M_k}\right\|\right) + 2 \left\| \mathbb{E}\left[ Q^TT\right]\right\|,\end{equation} and hence \begin{equation}\label{truegradatXkL2} \mathbb{E}\left[\left\|\nabla f(X_k)\right\|^2\right] \leq 4\sqrt{2}\left(\left\| \mathbb{E}\left[\Lambda^{-1}_{M_k}(X_k) \right\|^4\right]\right)^{1/2} \left( \epsilon_k^4 + 16\mathbb{E}\left[\left\| Q^TT\right\|^4\right]\right)^{1/2} + 8 \left\| \mathbb{E}\left[ Q^TT\right]\right\|^2.\end{equation} Also, some algebra yields that \begin{align}\label{nrdrineqs} \left| M_k^{-1}\sum_{j=1}^{M_k} Q_j^{(i)}Q_jX_k  \right|  \geq \left| M_k^{-1}\sum_{j=1}^{M_k} Q_j^{(i)}T_j  \right| - \epsilon_k \mbox{ and } \left| \mathbb{E}\left[Q^{(i)}Q\right]X_k  \right| \leq \left| \mathbb{E}\left[Q^{(i)}T \right] \right| + \epsilon_k. \end{align} Using~\eqref{nrdrineqs} in~\eqref{Lambda}, we have that the $i$-th diagonal element $\Lambda_{M_k}^{-1}(X_k)(i,i)$ of the matrix $\Lambda^{-1}_{M_k}(X_k)$ satisfies \begin{align}\label{lambdainvub} \mathbb{E}\left[ \left| \Lambda^{-1}(X_k)(i,i)\right|^4\right] & \leq 16\left( \left| \mathbb{E}\left[TQ^{(i)}\right]\right|^4 + \epsilon_k^4\right) \mathbb{E}\left[ \left( \left| M_k^{-1}\sum_{j=1}^{M_k} Q_j^{(i)}T_j  \right| - \epsilon_k \right)^{-4} \right] \nonumber \\ & < \infty,\end{align} where we have used the technical condition in~\eqref{techcondregression} to obtain the last inequality in~\eqref{lambdainvub}. Use~\eqref{lambdainvub} in~\eqref{truegradatXkL2} to conclude that the sufficient condition~\eqref{graditerbd} for $L_2$ convergence holds.$\blacksquare$ 
  
  It is interesting that for the bounded gradients assumption~\eqref{graditerbd} to hold, the further (somewhat stringent) technical condition in~\eqref{techcondregression} is needed even though $L_1$ convergence does not imply almost sure convergence in general. Also, while the treatment in Example~\ref{ex:regression} was for linear regression, similar conclusions for polynomial regression hold in a straightforward manner upon introducing higher powers of $x^*$ in~\eqref{regressionmodel}.

\subsection{Example (Karhounen-Lo\`{e}ve Error Class).}\label{klclass} The Karhunen-Lo\`{e}ve (KL) representation~\cite{1991ghaspa}, analogous to the Fourier representation of a deterministic function, is a series representation of a stochastic process as a linear combination of orthonormal functions having random coefficients that are pairwise uncorrelated. Models of physical processes based on KL are among the most fruitful applications of the theory of stochastic processes, and have been extensively used within statistical estimation and finite-element settings~\cite{1991ghaspa} over the past few decades. For this reason, such models are also natural contexts where stochastic optimization using a method such as RA may be of direct interest.

As a simple illustration, suppose the error process $\{\nabla F(x,Y) - \nabla f(x), x \in \mathcal{X}\}$ is a second-order process, that is, $\mathbb{E}\left[\|\nabla F(x,Y) - \nabla f(x)\|^2\right]< \infty$ for each $x \in \mathcal{X}.$ Then, $\{\nabla F(x,Y) - \nabla f(x), x \in \mathcal{X}\}$ admits a KL representation~\cite{2015ale,1991ghaspa}, that is, we can write \begin{equation}\label{KL} \nabla F(x,Y) - \nabla f(x) = \sum_{j=1}^{\infty} \sqrt{\lambda_j} \, Z_j(Y) \, \phi_j(x), \quad x \in \mathcal{X},\end{equation} where $\phi_j: \mathcal{X} \to \mathbb{R}^d, j \geq 1$ form an orthonormal basis (eigenfunctions) that spans $L^2(\mathcal{X})$, the associated eigen values $\lambda_j, j \geq 1$ and coefficients $\{Z_j, j \geq 1\}$ satisfy \begin{equation}\label{KLcoeffs} \mathbb{E}[Z_j] = 0; \quad \mbox{and} \quad \mathbb{E}[Z_iZ_j] = \delta_{ij}.\end{equation} The convergence in~\eqref{KL} is $L_2$, that is, $$e_n(x) := \nabla F(x,Y) - \nabla f(x) - \sum_{j=1}^{n} \sqrt{\lambda_j} \, Z_j(Y) \, \phi_j(x), \quad x \in \mathcal{X}$$ satisfies $\lim_{n \to \infty} \mathbb{E}\left[\|e_n(x)\|^2\right] = 0$ for each $x \in \mathcal{X}$. 

Suppose also that the following three conditions hold: \begin{enumerate} \item[(i)] $\exists b < \infty$ such that $\forall j \geq 1$, the eigenfunctions $\phi_j$ satisfy $$\sup_{x \in \mathcal{X}} \frac{\|\phi_j(x)\|}{\|\nabla f (x) \| +c_0} \leq b;$$ \item[(ii)] $\sum_{j=1}^{\infty} \lambda_j < \infty$; and \item[(iii)] there exist $\delta(\epsilon)>0$, $\sigma < \infty$ such that for each $j$: $$\mathbb{E}\left[ \left| \bar{Z}_{j,m}\right|^{2+\epsilon} \right] \leq \frac{\sigma}{m^{1+\delta(\epsilon)}}.$$ \end{enumerate} 

The assumption in (iii) holds under standard assumptions leading to the Bernstein inequality~\cite[pp. 33]{2018ver}, that is, when $\{Z_j, j \geq 1\}$ are independent, sub-exponential random variables; there exist extensions to the context where $\{Z_j, j \geq 1\}$ are sub-exponential but dependent --- see, for example,~\cite{2009dey} and~\cite[Section 5]{2018ver}.

Under (ii) and (iii), it can be shown that the sequence $\{\sum_{j=1}^n \sqrt{\lambda_j} Z_j(Y), n \geq 1\}$ is uniformly integrable and hence the equality in~\eqref{KL} holds pointwise (almost surely). Also, from~\eqref{KL}, we see that for each $x \in \mathcal{X}$, \begin{align}\label{KLfm} \|R_m(x)\| = \frac{\|\nabla f_m(x) - \nabla f(x)\|}{\|\nabla f(x)\| + c_0} &= \frac{1}{\|\nabla f (x)\| + c_0}\, \times  \left\|\frac{1}{m}\sum_{i=1}^m\sum_{j=1}^{\infty} \sqrt{\lambda_j}\, Z_j(Y_i) \, \phi_j(x)\right\|, \nonumber \\ & = \frac{1}{\|\nabla f (x)\| + c_0}\, \times \left\|\sum_{j=1}^{\infty} \sqrt{\lambda_j}\,\bar{Z}_{j,m}  \, \phi_j(x)\right\| \nonumber \\ & \leq \frac{b}{c_0}\, \sum_{j=1}^{\infty} \sqrt{\lambda_j}\,\left|\bar{Z}_{j,m} \right |,\end{align} where the last inequality in~\eqref{KLfm} follows from (i) above. Now notice that \begin{align}\label{bcasconvex} \sum_{m=1}^{\infty} P\left(\sum_{j=1}^{\infty} \sqrt{\lambda_j}\left | \bar{Z}_{j,m} \right| > \nu \right) & \leq \sum_{m=1}^{\infty} \frac{1}{\nu^{1+\epsilon}} \mathbb{E}\left[ \left(\sum_{j=1}^{\infty} \sqrt{\lambda_j} \left| \bar{Z}_{j,m} \right| \right)^{1+\epsilon} \right] \nonumber \\ & \leq \frac{1}{\nu^{1+\epsilon}} \left(\sum_{j=1}^{\infty} \sqrt{\lambda_j}\right)^{\epsilon}\sum_{m=1}^{\infty} \sum_{j=1}^{\infty} \mathbb{E}\left[ \sqrt{\lambda_j} \, \left | \bar{Z}_{j,m} \right |^{1+\epsilon}\right] \nonumber \\ & \leq \frac{1}{\nu^{1+\epsilon}}  \left(\sum_{j=1}^{\infty}\lambda_j\right)^{\frac{1+\epsilon}{2}} \left(\sum_{m=1}^{\infty} \mathbb{E}\left[ \left | \bar{Z}_{j,m} \right |^{2+2\epsilon}\right]\right)^{1/2} \nonumber \\ & < \infty, \end{align}
where the last inequality in~\eqref{bcasconvex} follows from (ii) and (iii) above. We now conclude from~\eqref{bcasconvex} and Borel-Cantelli's first lemma~\cite[pp. 59,60]{1995bil} that \begin{equation}\label{asconvex} \sum_{j=1}^{\infty} \sqrt{\lambda_j} |\bar{Z}_{j,m}| \, \as \, 0.\end{equation} Using~\eqref{asconvex} in~\eqref{KLfm}, we see that \begin{equation}\label{uniferrineq2}\sup_{x \in \mathcal{X}}\|R_m(x)\|  \as 0,\end{equation} implying that Theorem~\ref{thm:consistency} guarantees strong consistency if RA is implemented on this class with $M_k \as \infty$ and $\epsilon_k \as 0$. 

Using arguments similar to those in~\eqref{abserrgrowth}--\eqref{termcritagain} of Example I, we can also conclude that the key condition in~\eqref{graditerbd} holds.    
$\blacksquare$

\subsection{Example (Power Error Class).} \label{PECEx}
The key assumption of the complexity result in Theorem~\ref{thm:itercomp} is that the envelope of the weighted relative error class $w\mathcal{R} := \{R(x)( \| \nabla f(x)\|^{\alpha} + \delta_0), x \in \mathcal{X}\}$ exhibits the CLT-scaling~\eqref{emperrenv-clt}. The following simple example provides the most intuition.

Suppose the random error function $\nabla F(\cdot,Y) - \nabla f(\cdot)$ is such that for some $\gamma_0 \in (0,1)$, \begin{equation}\label{mulerrclass2} \nabla F(x,Y) - \nabla f(x)= A(Y) \|\nabla f(x)\|^{1-\gamma_0} + B(Y), \quad x \in \mathbb{R}^d\end{equation} where $A(Y)$ and $B(Y)$ are mean zero sub-exponential random variables.

Let's now demonstrate that the class characterized by~\eqref{mulerrclass2} (along with sub-exponential $A(Y),B(Y)$ assumption) satisfies the postulates of Theorem~\ref{thm:itercomp}, that is, its weighted relative error class has an empirical-mean envelope satisfying~\eqref{emperrenv-clt}. Choosing $\alpha < \gamma_0,$ and denoting $\bar{A}_m = m^{-1} \sum_{j=1}^m A(Y_j)$, $\bar{B}_m = m^{-1} \sum_{j=1}^m B(Y_j)$  we see that \begin{align}\label{mulclass3basic} \left\|R_m(x)\left(\| \nabla f(x) \|^{\alpha} + \delta_0\right)\right\| & \leq \left\| \frac{\bar{A}_m \| \nabla f(x) \|^{1-\gamma_0}}{\|\nabla f(x) \| + c_0} \left(\| \nabla f(x) \|^{\alpha} + \delta_0\right) \right\| + \left\| \frac{\bar{B}_m}{\|\nabla f(x) \| + c_0} \left(\| \nabla f(x) \|^{\alpha} + \delta_0\right) \right\| \nonumber \\ & \leq (1+\delta_0) (1 \vee \frac{1}{c_0}) (|\bar{A}_m| + |\bar{B}_m|).  \end{align} We see from~\eqref{mulclass3basic} that the empirical-mean envelope $H_m := \left\|R_m(x)\left(\| \nabla f(x) \|^{\alpha} + \delta_0\right)\right\|$ satisfies \begin{align}\label{exitercompmainineq} \mathbb{E}[m^{1/2\alpha}H_m^{1/\alpha}] & \leq  \left((1+\delta_0) (1 \vee \frac{1}{c_0})\right)^{1/\alpha}2^{1/\alpha}m^{1/2\alpha}\left(\mathbb{E}\left[|\bar{A}_m|^{1/\alpha}\right] + \mathbb{E}\left[|\bar{B}_m|^{1/\alpha}\right]\right) \nonumber \\
&\leq \left(2(1+\delta_0) (1 \vee \frac{1}{c_0})\right)^{1/\alpha}m^{1/2\alpha}\left(\mathbb{E}\left[|\bar{A}_m|^{1/\alpha}\right] + \mathbb{E}\left[|\bar{B}_m|^{1/\alpha}\right]\right). \end{align} Since $A(Y)$ is a sub-exponential random variable, and since $\bar{A}_m$ is the sample mean of iid random variables, there exist constants $c_A, \sigma_A \in (0,\infty)$ such that  \begin{align}\label{subexp} \mathbb{E}\left[\left|\bar{A}_m\right|^{1/\alpha}\right] &= \int_0^{\infty} P\left(\left| \sum_{j=1}^m A(Y_j) \right| > mt^{\alpha}\right) \nonumber \\ &\leq \int_0^{\sigma^2_Ac_A} 2\exp\{-\frac{1}{2}\frac{m}{\sigma^2_A}t^{2\alpha}\}\,dt +  \int_{\sigma^2_Ac_A}^{\infty} 2 \exp\{-\frac{1}{2}mc_At^{\alpha}\}\, dt \nonumber \\ &\leq \frac{1}{m^{1/2\alpha}}\left(\kappa_{1A} + \frac{\kappa_{2A}}{m^{1/2\alpha}}\right),\end{align} where $$\kappa_{1A} := \frac{2}{\alpha}\left(\frac{2}{c_A}\right)^{1/\alpha}\operatorname{\Gamma}(\frac{1}{\alpha}); \quad \kappa_{2A} := \frac{1}{\alpha}\left(\frac{2\sigma^2_A}{c_A}\right)^{1/2\alpha}\operatorname{\Gamma}(\frac{1}{2\alpha}),$$ and the second inequality in~\eqref{subexp} follows the standard tail bound on sub-exponential random variables~\cite{2018ver}. Similarly, since $B(Y)$ is a sub-exponential random variable, and since $\bar{B}_m$ is the sample mean of iid random variables, there exist constants $c_B, \sigma_B \in (0,\infty)$ such that  \begin{align}\label{subexpagain} \mathbb{E}\left[\left|\bar{B}_m\right|^{1/\alpha}\right] \leq \frac{1}{m^{1/2\alpha}}\left(\kappa_{1B} + \frac{\kappa_{2B}}{m^{1/2\alpha}}\right),\end{align} where $$\kappa_{1B} := \frac{2}{\alpha}\left(\frac{2}{c_B}\right)^{1/\alpha}\operatorname{\Gamma}(\frac{1}{\alpha}); \quad \kappa_{2B} := \frac{1}{\alpha}\left(\frac{2\sigma^2_B}{c_B}\right)^{1/2\alpha}\operatorname{\Gamma}(\frac{1}{2\alpha}).$$ Plugging~\eqref{subexp} and~\eqref{subexpagain} in~\eqref{exitercompmainineq}, we see that $\mathbb{E}[m^{1/2\alpha}H_m^{1/\alpha}]$ is finite if $\alpha < \gamma_0$ and we conclude that the power error class considered in this example satisfies the postulates of Theorem~\ref{thm:itercomp}. $\blacksquare$

\section{CENTRAL LIMIT THEOREM (CLT)}
We next treat the question of whether the sequence $\{\|\nabla f(X_k)\|, k \geq 1\}$ exhibits a CLT. Unlike the complexity results from Section~\ref{sec:itercomp} and Section~\ref{sec:oraccomp} which are non-asymptotic in the sense that they are valid for all $k \geq 1$, the CLT we present here holds only asymptotically for $k\gg 1$. CLTs say less than the complexity results about the behavior of $\{\|\nabla f(X_k)\|, k \geq 1\}$ for small $k$, but say more for large $k$ in the sense that they characterize the exact asymptotic distribution of $\|\nabla f(X_k)\|$. In particular, they precisely quantify the effect of the curvature of the objective function $f$ at the critical points on the rate of convergence. For this reason, CLTs are crucial building blocks for constructing confidence regions especially in the service of terminating an algorithm.    

Before we state the main result, we state an auxiliary result (without proof) that demonstrates the regular behavior of the scaled empirical error $\sqrt{M_k}(\nabla f_{M_k}(X_k) - \nabla f_{M_k}(X^*))$, where the random variable $X^*$ is such that $X_k \as X^* \in \mathcal{X}^*$. A proof of  Lemma~\ref{lem:aseq} follows along lines similar to that of~\cite[Theorem 6.3]{2006gee}. 
\begin{lemma}\label{lem:aseq} Suppose there exists a non-random non-increasing function $H: [0,1] \to \mathbb{R}^+$ such that \begin{equation}\label{entropybd} \lim_{A \to \infty} \limsup_{n \to \infty}P\left(\sup_{u >0} \frac{H_1(u,\mathcal{R},P_n)}{H(u)} > A\right)=0; \quad \int_0^1 H^{1/2}(u)\, du < \infty. \end{equation}\end{lemma} Define the $\delta$-ball around the set $\mathcal{X}^*$ of first-order critical points :  $$\mathcal{X}^*(\delta):= \bigcup_{x^* \in \mathcal{X}^*} \left\{x \in \mathbb{R}^d: \mathbb{E}\left[\|R(x,Y) - R(x^*,Y)\|\right] \leq \delta\right\}, \quad \delta>0.$$ \begin{enumerate} \item[(a)] There exists $\delta>0$ such that for any $x^* \in \mathcal{X}^*$, \begin{equation} \lim_{n \to \infty} P\left(\sup_{x \in \mathcal{X}^*(\delta) } \sqrt{n}\| \nabla f_n(x) - \nabla f_n(x^*)\| > \eta \right) < \eta.\end{equation}  \item[(b)] If the sequence of $\{X_k, k \geq 1\}$ of RA's iterates is such that $X_k \as X^*$ for some well-defined random variable $X^*$, then \begin{equation}\label{ascont} \sqrt{M_k} \| \nabla f_{M_k}(X_k) - \nabla f_{M_k} (X^*)\| \inD 0.\end{equation} \end{enumerate} 



We are now ready to state the main result of this section. In the result, we assume that the postulates of Theorem~\ref{thm:consistency} hold, implying that $\|\nabla f(X_k) \| \as 0$. To avoid contexts where $\|\nabla f(X_k) \| \as 0$ but the sequence $\{X_k, k \geq 1\}$ does not converge to anything, we also assume that there exists a random variable $X^*$ such that $X_k \as X^*$.

\begin{theorem}[Central Limit Theorem]
\label{thm:randcltXknonconvex}
Let the postulates of Theorem~\ref{thm:consistency}(b) hold along with the condition \begin{equation}\label{errtolcondb} \sqrt{M_k}\, \epsilon_k \inD 0. \tag{C.2(b)} \end{equation}  Suppose also that $X_k \as X^*
\in \mathcal{X}^*$ as $k \to \infty$. Denoting the covariance matrix $$\Sigma_x : = \int \nabla F(x,y) \nabla F(x,y)^{T}\, P(dy), \quad x \in \mathcal{X}^*.$$ Then we have as $k \to \infty$,
    \begin{equation}\label{randcltgrad}
        \sqrt{M_{k}} \, \nabla f(X_{k}) \inD N\left(0, \Sigma_{X^*}\right).\end{equation} If we further suppose that the objective function $f: \mathbb{R}^d \to \mathbb{R}$ is twice differentiable in $\mathcal{X}^*$ with invertible second-derivative matrix, and denoting $$V(X^*) : = \nabla^2 f(X^*)^{-1} \,\, \Sigma_{X^*} \,\, (\nabla^2 f(X^*)^{-1})^{T}.$$ Then we have that as $k \to \infty$, \begin{equation}\label{crlb} \sqrt{M_k} (X_k - X^*) \inD N(0,V(X^*)).\end{equation}
\end{theorem}
\proof{Proof.} Observe that 
\begin{align} \label{mainsplit}
    \sqrt{M_k} \nabla f(X_k) &= \sqrt{M_k}\left( \nabla f(X_k) - \nabla f_{M_k}(X_k) \right) + \sqrt{M_k}\nabla f_{M_k}(X_k) \nonumber \\ 
    &= \sqrt{M_k} \int \nabla F(X_k,Y) \, d (P - P_{M_k})  + \sqrt{M_k} \nabla f_{M_k}(X_k)  \nonumber \\
    &  = \sqrt{M_k} \int (\nabla F(X_k,Y) - \nabla F(X^*,Y)) \, d (P - P_{M_k}) \nonumber \\ & \hspace{1in} + \sqrt{M_k} \int \nabla F(X^*,Y) \, d (P - P_{M_k}) + \sqrt{M_k} \nabla f_{M_k}(X_k). 
\end{align} Lemma~\ref{lem:aseq} guarantees that  \begin{equation}\label{asequiseq}\sqrt{M_k} \int (\nabla F(X_k,Y) - \nabla F(X^*,Y)) \, d (P - P_{M_k}) \inP 0,\end{equation} and condition C.2(a) ensures that \begin{equation}\label{errtolcritagain}  \sqrt{M_k} \nabla f_{M_k}(X_k) \inD 0,\end{equation} We now plug~\eqref{asequiseq} and~\eqref{errtolcritagain} in~\eqref{mainsplit} to get \begin{equation}\label{randcltgradagain}
        \sqrt{M_{k}} \, \nabla f(X_{k}) \inD N\left(0, \Sigma(X^*)\right), \end{equation} and the first assertion of the theorem (appearing in~\eqref{randcltgrad}) holds.

To prove the second assertion (appearing in~\eqref{crlb}), we see that since $f$ is twice differentiable on $\mathcal{X}^*$ with a non-singular second-derivative $\nabla^2 f$, there exists a neighborhood $B(X^*,\epsilon)$ where \begin{equation}\label{secondderdef} \nabla f(X_k) = \nabla f(X_k) - \nabla f(X^*) = \nabla^2 f(X^*)(X_k - X^*) + o_p(\|X_k - X^*\|).\end{equation} On $B(X^*, \epsilon)$, we can then can write \begin{align}\label{emprocessscaled2} \sqrt{M_k} \nabla f(X_k) &= \sqrt{M_k} \nabla^2 f(X^*) (X_k - X^*) + o_p(\sqrt{M_k}\, \|X_k - X^*\|).\end{align} Now use~\eqref{randcltgradagain} in~\eqref{emprocessscaled2} to see that~\eqref{crlb} holds.$\blacksquare$
\endproof

The reader will observe that the random CLT in Theorem~\ref{thm:randcltXknonconvex} is different in nature than Theorem 12 in~\cite{2015kimpashen} (also see~\cite{2009shadenrus}) which is on the set of global minimizers, with the implication that its application to constructing confidence sets entails the practically difficult task of obtaining a global minimum estimator. Theorem~\ref{thm:randcltXknonconvex} is also different from the classical CLT in stochastic approximation, e.g., Theorem 1 in~\cite{1992poljud}, which in effect assumes strong convexity. Theorem~\ref{thm:randcltXknonconvex} should instead be seen as being analogous to Theorem 5.2 in~\cite{2011paskim} (see also~\cite[Ch. 10]{2003kusyin}).

\section{Numerical Experiments}\label{sec:numerical}
To investigate the empirical performance of RA, we ran a setting of RA alongside Adam or SGD (two highly popular stochastic gradient based methods) in three different contexts: least-squares, logistic regression, and a convolutional neural net (CNN), where the objective functions are strongly convex, convex, and nonconvex, respectively.

All algorithms were implemented in Python. The Tensorflow library was used for constructing the models and for computing gradients. Five different random runs were performed for each algorithm in each experimental setting. 
The same random seed initialization was set at the first run for each algorithm. 

To compare performance of the algorithms, we measure loss as a function of the total work performed by each algorithm. We display this in a single plot for each experiment: The y-axis of each plot displays $\log(f(x) - f^*)$, where $f(x)$ is the true loss function evaluated at $x$, and $f^*$ is the ``optimal'' loss value achieved by running a deterministic solver on the true loss for many iterations until convergence $(\|\nabla f(x)\|_2 \leq 10^{-8})$. The solid lines indicate the median loss, and the shaded region indicates the interquartile range of the losses from the five repeated runs. The loss is displayed as a function of the total number of individual gradients evaluated on the x-axis (referred to as ``Total Work"). For example, if RA takes $5$ inner iterations with a batch size of $100$, then this is measured as $500$ individual gradient evaluations. This is done so that RA and SGD/Adam are compared on similar ground, as RA is typically much more computationally intensive on a per-iteration basis.

\subsection{Models and Datasets}

For the least-squares experiment, we generated a dataset of 60,000 samples with 500 predictor variables. The response variable was set to be a linear function of the explanatory variables with i.i.d. standard Gaussian noise. The condition number of the covariance matrix (Hessian of the true objective function) was approximately $10^5$. Mean squared error of the residual error terms was used as the objective function.

The logistic regression experiment used the fashion MNIST dataset (consisting of 60,000 $28 \times 28$ greyscale images with 10 possible classes), and the CNN experiment used the LeNet CNN with the Cifar-10 dataset (60,000 $32 \times 32$ RBG images with 10 possible classes). For both experiments, cross-entropy loss was used with an $l_2$ penalty term of $\lambda = 1/n$, where $n$ is the number of observations in each dataset.

\subsection{Algorithms}

The implementation of RA was identical across each experiment. We used a sample increase rate of \rv{$c_1:=1.05$ ($M_{k+1} = c_1 M_k, \forall k$)}, with a starting batch size of $M_0:=32$. L-BFGS with backtracking line search \cite[Algorithm 7.5]{nocedalbook} was used as the deterministic solver for the inner iterations. Updates to BFGS matrices are skipped if the \emph{curvature condition} is not satisfied (see \cite{BollapragadaICML18} ). Notably, the stored gradient and iterate differences used in L-BFGS were carried across outer iterations instead of restarting L-BFGS for every new sample-path problem. 

The threshold $\epsilon_k$ for terminating L-BFGS during the $k^{th}$ iteration was set to be $\hat{\sigma}_k^2/\sqrt{M_k}$, where $\hat{\sigma}_k^2$ is defined to be the trace of the sample gradient covariance matrix of $\nabla F_{S_\sigma}(X_k) $ for some subset of RA's current sample of size $S_\sigma := \min \{128, M_k\}$.

For implementations of Adam and SGD, the best two step sizes (in terms of final loss) were chosen from the set of $\{ 10^{j} \},j\in \Z$. Step sizes were held constant throughout each run. For Adam's additional two hyper-parameters, the default Tensorflow values were used. Batch sizes for SGD and Adam were set at the default 32.

\subsection{Least-Squares}

\begin{figure}[h]
\centering
    \includegraphics[scale=.7]{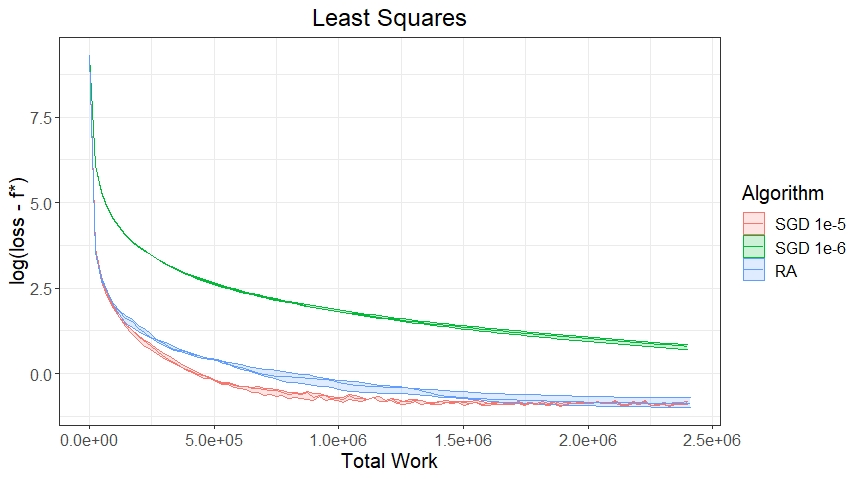}
    \caption{\textbf{Least Squares}}
    \label{fig:linreg}
\end{figure}

Results for the least squares experiment displayed in Figure \ref{fig:linreg}. The optimally-tuned SGD (step size of $10^{-5}$) performs similarly to RA, with both achieving comparable loss value after 40 epochs. The loss obtained by the second best step size for SGD is substantially higher than that of RA and the best tuned SGD.

\subsection{Multiple Logistic Regression}

\begin{figure}[h]
    \includegraphics[scale=.7]{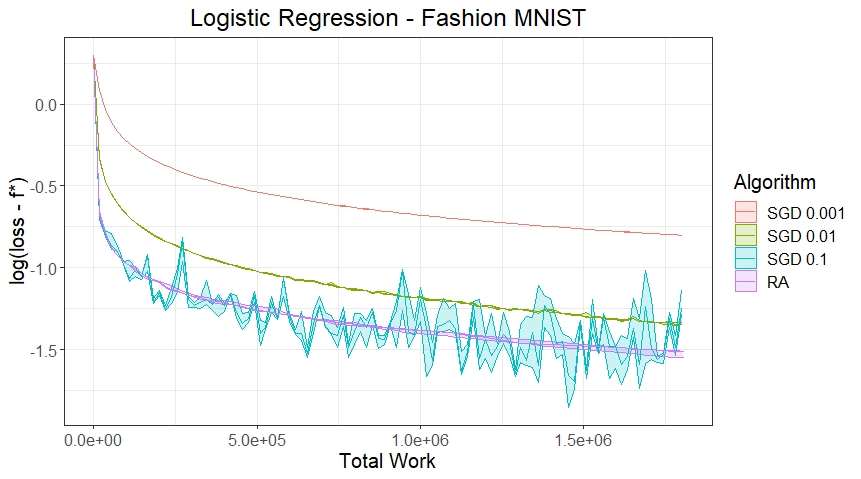}
    \caption{\textbf{Logistic Regression on Fashion MNIST}}
    \label{fig:logreg}
\end{figure}

Results for the logistic regression experiment are displayed in Figure \ref{fig:logreg}. The optimally-tuned SGD (step size of $10^{-2}$) achieves similar loss to RA, but has considerably more variability with respect to the loss values. Using a smaller step size with SGD smooths out the algorithm's loss trajectory, but the descent is much slower as a result.

\subsection{LeNet on Cifar-10}

\begin{figure}[h]
    \includegraphics[scale=.7]{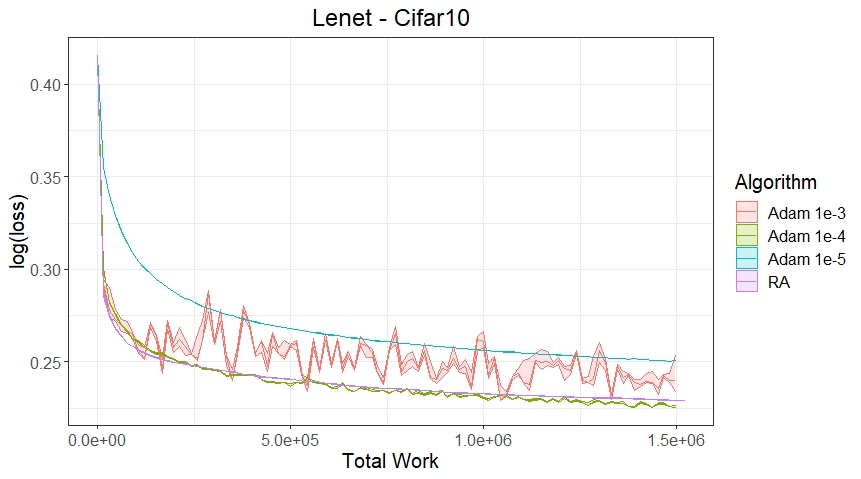}
    \caption{\textbf{LeNet on Cifar-10}}
    \label{fig:cnn1}
\end{figure}

Results for the CNN experiment  are displayed in Figure \ref{fig:cnn1}. The optimally-tuned Adam (step size of $10^{-4}$) achieves comparable loss to RA, but has slightly more variability with respect to the loss values. Step size parameters of $10^{-3}$ and $10^{-5}$ show deteriorating performance (in terms of variability and loss, respectively) for Adam.

\subsection{Discussion of Numerical Results}

\subsubsection{Hyper-parameter Tuning}

The results from the numerical experiments illustrate some distinct features of RA compared to Adam or SGD. One clear advantage of RA is the lack of need for hyper-parameter tuning, due to its incorporation of a deterministic solver on the sample-path problems. Adam's and SGD's performance, on the other hand, is highly dependent upon the user selecting an appropriate step size (in addition to any other parameters). This phenomenon is demonstrated across all experiments, where only one choice of step size in SGD/Adam generates a loss trajectory comparable to RA. Ideally, solving an optimization problem would not require the user to have to hand-tune parameters for every individual problem. Thus, it may be argued that RA dispenses with tuning in a stochastic optimization context, demonstrating robust performance across several settings.

\subsubsection{Asymptotic Consistency}
From a theoretical perspective, RA guarantees several types of convergence under fairly standard assumptions on the objective function and oracle. On the other hand, SGD and Adam converge (in expectation) only if their step-sizes approach 0 at an appropriate rate. However, this often results in much poorer finite-time performance, which is likely why Tensorflow's and PyTorch's implementations of SGD and Adam use constant step-sizes by default (although more advanced techniques, such as cycling the step-size or reducing it adaptively, are available). RA seems to bridge this gap by provably being able to converge while at the same time showing robust finite-time performance without any need for parameter tuning.

\section{Concluding Remarks} A widely practiced method of solving stochastic optimization problems is to execute a ``plug-in'' second-order solver such as L-BFGS line search or Newton-CG to a specified tolerance on a sample-path approximation of the objective function. Retrospective Approximation (RA) makes this idea rigorous: generate a sequence of approximate sample-path problems and solve (using an existing deterministic second-order solver) to within a sequence of decreasing error tolerances. Another way to think about RA is as a stochastic second-order method that allows for multiple optimization steps between sample updates of the gradient. To ensure optimal convergence rates, or in loose terms avoid ``oversolving'' of sample-path problems, RA stipulates a certain relationship between the error tolerances to which the sample-path problems are solved and the sample sizes using which the sample-path functions are generated. The optimal theoretical properties of RA hold when the relative-error class associated with the stochastic oracle in use satisfies conditions that arguably hold widely, and in particular are much weaker than what is typically imposed within the SAA literature.   

\section*{Acknowledgments.} Raghu Pasupathy gratefully acknowledges the Office of Naval Research for support provided by the grants N000141712295 and 13000991, and Prof. Susan Hunter at Purdue Industrial Engineering for insightful discussions. Raghu Bollapragada gratefully acknowledges the support of Lawrence Livermore National Laboratory and the National Science Foundation for support provided by the grant NSF DMS 2324643.

\bibliographystyle{informs2014} 
\bibliography{stochasticoptimization.bib,BigReferences.bib}

\newpage
\begin{center} \large{SUPPLEMENTARY MATERIAL}
\end{center}
\begin{APPENDICES}

\section{Proof of Theorem~\ref{thm:itercomp}}\label{ProofAltThm}

    Recall that \begin{equation}\label{redef}R_{M_k}(X_k) := \frac{\nabla f_{M_k}(X_k) - \nabla f(X_k)}{\|\nabla f(X_k)\| + c_0}.\end{equation} Taking norm and applying triangle inequality to~\eqref{redef}, we can write \begin{equation} \label{initsplit} \| \nabla f (X_k) \| \leq \| \nabla f_{M_k}(X_k)\| + c_0\|R_{M_k}(X_k) \| + \|R_{M_k}(X_k)\| \| \nabla f(X_k) \|.\end{equation} Take expectation on both sides of~\eqref{initsplit} and noticing that $\|\nabla f_{M_k}(X_k)\| \leq \epsilon_k$, we can write \begin{align}\label{secondsplit} \mathbb{E}\left[ \| \nabla f (X_k) \| \right] &\leq \epsilon_k + c_0 \mathbb{E}\left[\|R_{M_k}(X_k)\| \frac{\|\nabla f(X_k) \|^{\alpha} + \delta_0}{\|\nabla f(X_k) \|^{\alpha} + \delta_0} \right]  + \mathbb{E}\left[\|R_{M_k}(X_k)\| \|\nabla f(X_k) \|\frac{\|\nabla f(X_k) \|^{\alpha} + \delta_0}{\|\nabla f(X_k) \|^{\alpha} + \delta_0} \right] \nonumber \\ & \leq \epsilon_k + \frac{c_0}{\delta_0} \mathbb{E}\left[\|R_{M_k}(X_k)\| \left(\|\nabla f(X_k) \|^{\alpha} + \delta_0\right) \right]  \nonumber \\ & \hspace{2in} + \mathbb{E}\left[\underbrace{\|R_{M_k}(X_k)\| \left(\|\nabla f(X_k) \|^{\alpha} + \delta_0\right)}_{A} \underbrace{\| \nabla f(X_k) \|^{1-\alpha}}_{B} \right]. \end{align}  Recall from Young's inequality~\cite[pp. 13]{1989kre} that \begin{equation}\label{youngsineq} AB \leq \frac{\gamma^p A^p}{p} + \frac{B^q}{\gamma^q q}; \quad \frac{1}{p} + \frac{1}{q} = 1; \quad p,q \geq 1; \quad \gamma >0.\end{equation} Apply ~\eqref{youngsineq} to the product of the quantities designated $A$ and $B$ in~\eqref{secondsplit}, choose $$p: = \frac{1}{\alpha}; \quad q: = \frac{1}{1-\alpha}; \quad \gamma := 2(1-\alpha)^{1-\alpha},$$ and take expectation to get \begin{align}\label{youngsapp} \MoveEqLeft \mathbb{E}\left[\|R_{M_k}(X_k)\| \left(\|\nabla f(X_k) \|^{\alpha} + \delta_0\right) \, \| \nabla f(X_k) \|^{1-\alpha} \right]  \nonumber \\ & \leq \frac{\alpha(2(1-\alpha))^{\frac{1}{\alpha}}}{1-\alpha} \mathbb{E}\left[\left(\|R_{M_k}(X_k)\| \left(\|\nabla f(X_k) \|^{\alpha} + \delta_0\right)\right)^{1/\alpha}  \right] + \frac{1}{2} \mathbb{E}\left[\| \nabla f(X_k)\| \right].  \end{align} Plugging~\eqref{youngsapp} in~\eqref{secondsplit}, we get \begin{align}\label{keysplit} \frac{1}{2}\mathbb{E}\left[ \| \nabla f (X_k) \| \right] &\leq  \epsilon_k + \frac{c_0}{\delta_0} \mathbb{E}\left[\|R_{M_k}(X_k)\| \left(\|\nabla f(X_k) \|^{\alpha} + \delta_0\right) \right]  \nonumber \\ & \hspace{1in} + \frac{\alpha(2(1-\alpha))^{\frac{1}{\alpha}}}{1-\alpha} \mathbb{E}\left[\left(\|R_{M_k}(X_k)\| \left(\|\nabla f(X_k) \|^{\alpha} + \delta_0\right)\right)^{1/\alpha}  \right] \nonumber \\ & \leq \epsilon_k + \frac{c_0}{\delta_0} \mathbb{E}\left[\mathbb{E}\left[\sup_{x \in \mathbb{R}^d}\|R_{M_k}(x)\| \left(\|\nabla f(x) \|^{\alpha} + \delta_0\right) \, \vert \, \mathcal{F}_{k-1} \right]\right] \nonumber \\ & \hspace{0.5in} + \frac{\alpha(2(1-\alpha))^{\frac{1}{\alpha}}}{1-\alpha} \mathbb{E}\left[\mathbb{E}\left[\left(\sup_{x \in \mathbb{R}^d}\|R_{M_k}(x)\| \left(\|\nabla f(x) \|^{\alpha} + \delta_0\right)\right)^{1/\alpha} \, \vert \, \mathcal{F}_{k-1} \right]\right] \nonumber \\ & = \epsilon_k + \frac{c_0}{\delta_0} \mathbb{E}\left[\mathbb{E}\left[H_{M_k} \, \vert \, \mathcal{F}_{k-1} \right]\right] + \frac{\alpha(2(1-\alpha))^{\frac{1}{\alpha}}}{1-\alpha} \mathbb{E}\left[\mathbb{E}\left[H_{M_k}^{1/\alpha} \, \vert \, \mathcal{F}_{k-1} \right]\right].  \end{align} Due to the assumptionin~\eqref{emperrenv-clt}, and since $\alpha \in (0,1)$, we know that there exists $\sigma_0 < \infty$ such that \begin{equation}\label{VCineq} \mathbb{E}\left[\left(\sqrt{M_k} \, H_{M_k}\right)^n \, \vert \, \mathcal{F}_{k-1} \right] \leq  \sigma_0^n \mbox{ for } n \leq \frac{1}{\alpha}. \end{equation} Plugging~\eqref{VCineq} in~\eqref{keysplit} and since $\epsilon_k := C_{2,k}/M_k^p$ where  $\mathcal{F}_{k-1} \ni C_{2,k} \leq c_2$, we get \begin{align}\label{keysplitfinal} \mathbb{E}\left[ \| \nabla f (X_k) \| \, \vert \, \mathcal{F}_{k-1} \right] & \leq  \frac{2c_2}{M_k^p} + \frac{2c_0\sigma_0/\delta_0}{\sqrt{M_k}} +  \left(\frac{\kappa_2}{\sqrt{M_k}}\right)^{1/\alpha}, \end{align} proving~\eqref{geniter}. Now plug the choice of sample size sequence given in~\eqref{sampsizeseqchoice} in~\eqref{keysplitfinal} and $p = 1/2$ to see that the iteration complexity~\eqref{itercomp} holds as well.

\section{Key Results from the Theory of Empirical Processes}
In this section, we provide some relevant results that relate the 
\eqref{ulln}, or \eqref{ass:RULLN} to the underlying geometric or regularity 
property of the functions generated by the oracle called upon in our algorithm. 
Note that the crucial property is the concept
of ``uniform'' in the Law of Large Numbers within a class of functions. 
Such a theory has been developed
quite extensively in the realm of empirical processes.
We refer to \cite{1996vanwel, 2016ginnic, 2006gee} for comprehensive treatments.

Recall that by classical Law of Large Numbers, given any probability 
measure $P$ and function $F$ with $P|F| < \infty$, we have that
$$
P_nF := \frac{1}{n}\sum_{i=1}^n F(Y_i)
\longrightarrow_{n\rightarrow\infty} PF:=\mathbb{E}F(Y),
\quad \text{$P$ a.s.}
$$
where $Y_i$ are IID samples with distribution given by $P$. The notation
$P_n$ refers to the empirical measure with respect to the samples. We would 
like the 
above limit to hold uniformly for a class $\mathcal{F}$ of functions.
This is captured by the following uniform norm:
$$
\|P_n-P\|_{\mathcal{F}} := \sup_{F\in\mathcal{F}}\left|P_nF - PF\right|
$$
Then the Uniform Law of Large Numbers (ULLN) with respect to the class $\mathcal{F}$
is stated as:
$$
\lim_{n\rightarrow0}\|P_n-P\|_{\mathcal{F}} = 0,
\quad \text{$P$ a.s.}
$$
If the above is true, then $\mathcal{F}$ is called a \emph{Glivenko-Cantelli} 
class.

Translating the above framework to our current RA setting, 
$F$ and $\mathcal{F}$ are specifically given by
$$
\mathcal{F} = \left\{
F_x:=F(x,\cdot): x\in \mathbb{R}^d\,\,\,\text{(or $\mathbb{R}^d$)}\right\}
$$
so that the functions are indexed by $x\in\mathbb{R}^d$. We will further
consider the relative error class \eqref{relerror} repeated here for
convenience:
$$
\mathcal{R} : = \left\{R_x:= R(x,\cdot),
\quad x\in\mathbb{R}^d\,\,\,\text{(or $\mathbb{R}^d$)}
\right\}\quad\text{where}\quad
R(x,Y) := \frac{\nabla F(x,Y) - \nabla f(x)}{\| \nabla f(x) \| + c_0}.
$$
In the above cases, we have
\begin{equation*}
\left\{
\begin{array}{cl}
P_nF_x &= \frac{1}{n}\sum_{i=1}^n F(x,Y_i),\\
PF_x   &= \mathbb{E}F(x,Y)\,\,\,(:= f(x)),\\
\|P_n - P\|_{\mathcal{F}} 
& = \sup_{x\in\mathbb{R}^d}\left| P_nF_x - f(x)\right|
\end{array}
\right.
\quad\text{and}\quad
\left\{
\begin{array}{cl}
P_nR_x &= \frac{\frac{1}{n}\sum_{i=1}^n F(x,Y_i) - \nabla f(x)}
{\| \nabla f(x) \| + c_0},\\
PR_x   &= \mathbb{E}R(x,Y)\,\,\,(=0),\\
\|P_n - P\|_{\mathcal{R}} 
& = \sup_{x\in\mathbb{R}^d}\left| P_nR_x\right|
\end{array}
\right.
\end{equation*}
Hence, the Uniform Law of Large Numbers is exactly given by 
$$
\lim_{n\rightarrow\infty}\sup_{x\in\mathbb{R}^d}\left|
\frac{1}{n}\sum_{i=1}^n \nabla F(x,Y_i) - \nabla f(x)
\right|=0,
\quad\text{and}\quad
\lim_{n\rightarrow\infty}\sup_{x\in\mathbb{R}^d}\left|
\frac{\frac{1}{n}\sum_{i=1}^n \nabla F(x,Y_i) - \nabla f(x)}
{\| \nabla f(x) \| + c_0}
\right|
=0.
$$
But in order to keep this Appendix in a general setting, we will follow the 
usual notation in empirical processes, using $F$ and $\mathcal{F}$ to denote a 
particular function and the function class.

As mentioned before, the essential criterior for the validity of ULLN
is given by the entropy of $\mathcal{F}$ which measures the ``(ir)regularity''
or ``richness'' of the function class. This is given by the following 
definition.
\begin{definition}[Empirical Entropy]
A $\delta$-net for $\mathcal{F}$ is a set $\{F_1,F_2, \ldots,F_n\}$ of functions
 in $\mathcal{F}$ such that any $F$ can be approximated to within $\delta$ in 
$L_p(P)$ norm ($p\geq 1$) using one of the $F_j$'s.
That is, for any $F\in\mathcal{F}$, there is some 
$j \in \{1,2,\ldots,n\}$ such that, 
$$
\left(\mathbb{E}\bigg[ \left\| F(Y) - F_j(Y)\right\|^p\bigg]\right)^{1/p} 
= \left(\int\|F(Y) - F_j(Y)\|^p\, \diff P\right)^{1/p} \leq \delta.
$$ 

The \emph{$\delta$-covering number} $N_p(\delta,\mathcal{F},P)$ of the class 
$\mathcal{F}$ is then the ``smallest'' $\delta$-net,
i.e.,  $$N_p(\delta,\mathcal{F},P) := \inf\bigg\{n: \mbox{there exists a } \delta\mbox{-net of } \mathcal{F} \mbox{ with cardinality } n\bigg\}. $$ 

The \emph{entropy} of the class $\mathcal{F}$ is the logarithm of its $\delta$-covering number: 
$$ H_p(\delta,\mathcal{F},P) = \log N_p(\delta,\mathcal{F},P).$$
\end{definition}

In actual application, the empirical version of the above is often used. 
More precisely, the measure $P$ is replaced by its empirical measure $P_n$, 
leading to the following \emph{empirical entropy}: 
$$
H_p(\delta,\mathcal{F},P_n).
$$
Note that the above is a \emph{random} quantity. We further recall the 
definition here of the envelope $G$ given in \eqref{defn:env}:
\begin{equation}
G = G(Y) := \sup \left\{\|F(Y)\|: F \in \mathcal{F}\right\}.
\tag{\text{\ref{defn:env}}}
\end{equation}
It is interesting to note that by just assuming $\mathbb{E}G < \infty$, we already
have that $\|P_n - P\|_{\mathcal{R}}$ converges to some \emph{finite limit}
$P$ a.s. and in $L^1$ \cite[Prop. 3.7.8]{2016ginnic}. 

To show that the limit
is actually zero, a growth restriction on $H_p$ is crucial. More precisely,
let 
$$
\mathcal{F}_B = \left\{
F\mathbb{I}_{F\leq B}: F\in\mathcal{F}
\right\},\quad\text{where}\,\,\,
F\mathbb{I}_{F\leq B}=\left\{
\begin{array}{ll}
F(x, Y), & \|F(x,Y)\| \leq M,\\
0, & \|F(x,Y)\| > M.
\end{array}
\right.
$$
Then by \cite[Thm. 3.7.14]{2016ginnic}, under the assumption that 
$\mathcal{F}$ is $L^1$-bounded, i.e.
$\sup_{F\in\mathcal{F}}\mathbb{E}|F| < \infty$,
the fact that $\mathcal{F}$ being Glivenko-Cantelli
is equivalent to the following statement: 
\begin{center}
$\mathbb{E}G < \infty$, and for all $B < \infty$ and $\delta > 0$, 
and some $p \in (0,\infty]$, it holds that
\begin{equation}\label{ass:sublinentgrowth} 
\frac{1}{n}H_p(\delta,\mathcal{F},P_{n}) \inP 0 
\mbox{ as } n \to \infty. \tag{Sub-Ent}
\end{equation}
\end{center}
By \cite[Corollary 3.7.15]{2016ginnic}, the above is also equivalent to
\begin{equation}\label{ass:sublinengrowth2}
\text{$\mathbb{E}G < \infty$ and}\,\,\,
\frac{1}{n}\mathbb{E}\sqrt{H_2(\delta, \mathcal{F}, P_n)} \longrightarrow 0.
\tag{Sub-Ent$_2$}
\end{equation}
See also \cite[Theorem 2.4.3]{1996vanwel} for a similar statement.

We make the following remarks pertaining to the analysis of RA.
\begin{enumerate}
\item[(a)]
An implication of such equivalence is that the main condition~\eqref{ass:RULLN} in Theorem~\ref{thm:consistency} can be replaced by~\eqref{ass:sublinentgrowth} without changing the assertions of the theorem. Even though~\eqref{ass:sublinentgrowth} is one of the key ``workhorse'' conditions in empirical process theory, we have chosen to present our results with~\eqref{ass:RULLN} since, as we shall demonstrate in Section~\ref{sec:examples}, the condition in~\eqref{ass:RULLN} appears to be much easier to check in optimization contexts. 

\item[(b)] Following the proofs of
\cite[Thm. 3.7.14, Corollary 3.7.15]{2016ginnic} and
\cite[Theorem 2.4.3]{1996vanwel}, we can extend the results to 
second or higher order moments and with respect to a filtration.
In particular, using the notations of Theorems \ref{thm:consistency} and 
\ref{thm:l1cons}, under the assumption 
$\mathbb{E}\left[G\right] < \infty,$ the sublinear entropy growth 
condition~\eqref{ass:sublinentgrowth} can be used to show that
$$
\sup_x \|R_{M_k}(x)\| \as 0
\quad\text{and}\quad
\mathbb{E}\left[\sup_x \|R_{M_k}(x)\| \, \vert \, \mathcal{F}_{k-1}\right] 
\to 0
\quad\text{if}\quad
M_k \as \infty.
$$ 
Furthermore, if $ \mathbb{E}\left[G^2\right] < \infty,$ 
then we can have
$$\mathbb{E}\left[\sup_x \|R_{M_k}(x)\|^2 \vert \, \mathcal{F}_{k-1}\right] 
\as 0
\quad\text{if}\quad
M_k \as \infty.
$$
This is something we will use in Theorem~\ref{thm:l1cons} when we demonstrate 
the $L_1$-consistency.
\end{enumerate}

The next question of theoretical and practical importance is the \emph{rate}
of convergence in \eqref{ulln}. By the classical Central Limit Theorem (CLT), 
we expect $\|P_n-P\|_{\mathcal{F}} \lesssim O(\frac{1}{\sqrt{n}})$. Again, the
uniformity in $\mathcal{F}$ is crucial in our application. This is achieved
by means of the following quantity \cite[p. 186]{2016ginnic}.
\begin{definition}[Koltchinskii-Pollard entropy]
Given a class of function $\mathcal{R}$ and its envelope function $G$,
\begin{equation}
J(\mathcal{F}, G, \delta)
= \int_0^\delta
\sup_{Q}\sqrt{
H_2\left(\tau\|G\|_{L^2(Q)}, \mathcal{F}, Q\right)
}
\,d\tau
\end{equation}
where the $\sup$ is taken over all finite discrete probabilities $Q$ with
$\|G\|_{L^2(Q)} > 0$.
\end{definition}

With the above, we have the following theorem which exactly gives the
type of CLT scaling we need \cite[Theorem 3.5.4]{2016ginnic}.
\begin{theorem} The following estimate holds for all $n\in\mathbb{N}$:
\begin{equation}
\mathbb{E}\left[
\sqrt{n}\|P_n-P\|_{\mathcal{F}}
\right]
\leq
\max\left\{
A_1\|G\|_{L^2(P)}J(\mathcal{F}, G, \delta), 
\frac{A_2\|U\|_{L^2(P)}J^2(\mathcal{F}, G, \delta)}{\sqrt{n}\delta^2}
\right\}
\end{equation}
where 
$$
U = \max_{1\leq i \leq n} G(Y_i),
\quad
\sigma^2 = \sup_{F \in\mathcal{F}}PF^2,
\quad
\delta = \frac{\sigma}{\|G\|_{L^2(P)}},
$$
and $A_1$ and $A_2$ are some explicit universal constants.
\end{theorem}

\begin{remark}
Note that all of the above results about the limiting behavior of 
empirical measures make use of the entropy function $H_p$. 
Readers might be aware of another concept which can provide explicit bounds on 
$H_p$. It is given by means of VC dimension, named after 
Vladimir Vapnik and Alexey Chervonenkis. VC dimension measures the combinatorial 
richness of the class of functions in question. Without getting too much into detail, we will simply state the following definition which illustrates the 
type of bounds we can have on $H_p$.
\begin{definition}\cite[Definition 3.6.10]{2016ginnic} A class of measurable 
functions $\mathcal{F}$ is of \emph{VC} type with respect to a envelope function
$G$ if there exist finite constants $A$ and $\nu$ such that for all probability
meaures $Q$ and $0 < \epsilon < 1$, it holds that
$$
N_2(\epsilon\|G\|_{L^2(Q)}, \mathcal{F}, Q) 
\leq \left(\frac{A}{\epsilon}\right)^\nu.
$$
\end{definition}
The above estimate clealy provides sufficient bound for the entropy function
that will be workable for all our results.
\end{remark}

The following concentration inequality by Talagrand, though not needed in the current work, plays a fundamental role in 
the theory of empirical processes. 

\begin{theorem}\cite[Theorem 1.4]{talagrand1996new}
Using the previous notations, in particular, 
let $Y_i, \ldots, Y_n$ be IID variables distributed according to $P$. 
Let $\mathcal{F}$ be any class of functions that is uniformly bounded by a 
constant $U > 0$. Then for all $t > 0$, 
\begin{equation}\label{tala.conc.ineq}
P\left\{
\left|
\left\|\sum_{i=1}^n F(Y_i)\right\|_{\mathcal{F}}-
\mathbb{E}\left\|\sum_{i=1}^n F(Y_i)\right\|_{\mathcal{F}}
\right| \geq t
\right\}
\leq 
K \exp\left\{
-\frac{1}{K}\frac{t}{U}\log\left(
1 + \frac{tU}{V}
\right)
\right\}
\end{equation}
where $K$ is a universal constant and $V$ is any number satisfying
$$
V \geq \mathbb{E}\sup_{F\in\mathcal{F}}
\sum_{i=1}^n F^2(X_i)
$$
\end{theorem}
In order to see the explicit scaling on $n$, we can use $V = nU^2$. Upon 
replacing $t$ by $\frac{t}{n}$, then \eqref{tala.conc.ineq} becomes
\begin{equation}\label{tala.conc.ineq.n}
P\left\{
\left|
\left\|\frac{1}{n}\sum_{i=1}^n F(Y_i)\right\|_{\mathcal{F}}-
\mathbb{E}\left\|\frac{1}{n}\sum_{i=1}^n F(Y_i)\right\|_{\mathcal{F}}
\right| \geq t
\right\}
\leq
K \exp\left\{
-\frac{1}{K}\frac{nt}{U}\log\left(
1 + \frac{t}{U}
\right)
\right\}.
\end{equation}
This is the form used in the proof of Theorem \ref{thm:l1cons}.

\end{APPENDICES}

\end{document}